\documentclass[aps,prx,superscriptaddress,twocolumn]{revtex4-2}
\usepackage{hyperref}
\usepackage{graphics}
\usepackage{graphicx} 
\usepackage{rotating}
\usepackage{dcolumn}
\usepackage{longtable}
\usepackage{multirow}
\usepackage{amsmath}
\usepackage{amsthm}
\usepackage{amssymb}
\usepackage{gensymb}
\usepackage{centernot}
\usepackage{xspace}
\usepackage{textcase}
\usepackage[margin=1in]{geometry}
\usepackage{natbib}
\usepackage{float}
\usepackage{tikz}
\usetikzlibrary{trees}
\usetikzlibrary{arrows}
\usepackage{bbm}
\usepackage{url}
\usepackage[ruled,vlined]{algorithm2e}
\SetKwInOut{Input}{input}\SetKwInOut{Output}{output}

\newtheorem{theorem}{Theorem}[section]
\newtheorem*{theorem*}{Theorem}
\newtheorem{proposition}{Proposition}[section]
\newtheorem*{proposition*}{Proposition}

\newtheorem{lemma}[theorem]{Lemma}
\newtheorem*{lemma*}{Lemma}
\theoremstyle{definition}
\newtheorem{definition}{Definition}[section]
\newtheorem{remark}{Remark}
\newtheorem{example}{Example}

\newcommand{\state}{\omega}
\newcommand{\statespace}{\Omega}
\newcommand{\Reals}{\mathbb{R}}
\newcommand{\Koop}{\mathcal{K}}
\newcommand{\KoopMat}{K}

\newcommand{\M}{\Omega}
\newcommand{\F}{\mathcal{F}}

\newcommand{\rossler}{R\"{o}ssler\xspace}
\newcommand{\causes}{\mathbf{\longrightarrow}}
\newcommand{\ncauses}{\mathbf{\centernot\longrightarrow}}
\newcommand{\ddcauses}{\xrightarrow[\KoopMat^t]{}}

\begin{document}

\title{Causal Discovery in Nonlinear Dynamical Systems using Koopman Operators}

\author{Adam Rupe}
\affiliation{Pacific Northwest National Laboratory}
\affiliation{EnviTrace LLC}
\author{Derek DeSantis}
\affiliation{Los Alamos National Laboratory}
\author{Craig Bakker}
\affiliation{Pacific Northwest National Laboratory}
\author{Parvathi Kooloth}
\affiliation{Pacific Northwest National Laboratory}
\author{Jian Lu}
\affiliation{Pacific Northwest National Laboratory}
\date{\today}

\begin{abstract}
    We present a theory of causality in dynamical systems using Koopman operators. Our theory is grounded on a rigorous definition of causal mechanism in dynamical systems given in terms of flow maps. In the Koopman framework, we prove that causal mechanisms manifest as particular flows of observables between function subspaces. While the flow map definition is a clear generalization of the standard definition of causal mechanism given in the structural causal model framework, the flow maps are complicated objects that are not tractable to work with in practice. By contrast, the equivalent Koopman definition lends itself to a straightforward data-driven algorithm that can quantify multivariate causal relations in high-dimensional nonlinear dynamical systems. The coupled \rossler system provides examples and demonstrations throughout our exposition. We also demonstrate the utility of our data-driven Koopman causality measure by identifying causal flow in the Lorenz 96 system. We show that the causal flow identified by our data-driven algorithm agrees with the information flow identified through a perturbation propagation experiment. Our work provides new theoretical insights into causality for nonlinear dynamical systems, as well as a new toolkit for data-driven causal analysis.
\end{abstract}

\maketitle


\section{Introduction}

For many complex dynamical systems, there is a disconnect between the governing equations and system behaviors due to a web of strong nonlinear dependencies \cite{rupe24a}. Simulations that numerically solve the governing equations are just a starting point for understanding the complex behaviors that emerge. Causal discovery is a potent form of analysis to identify and quantify the detailed processes that give rise to these behaviors \cite{rung19a,camp23a}. 
The standard statistical framework for causal inference \cite{pearl09a}, however, is insufficient for identifying causal relations in nonlinear dynamical systems. 
While tools specifically designed for causal discovery in dynamical systems have shown promise, some notable gaps remain in their theoretical foundations and algorithmic implementation. In this work we introduce a novel framework for causal discovery in nonlinear dynamical systems using Koopman operators to address these shortcomings. 

Essentially all of the statistical approaches to causal inference assume that causal flow is unidirectional, represented by an acyclic causal graph, and many assume linear causal relationships. These assumptions do not hold for nonlinear dynamical systems. 
The unidirectional assumption cannot account for feedbacks, which are ubiquitous in complex nonlinear systems. In climate science, for example, a central question is how the Earth system responds to forcing~\cite{ghil20a}. One can simulate various input forcings and record the resulting responses, but it is not necessarily straightforward to follow all of the processes that lead from the input forcing to the output response. This is particularly the case when there are feedbacks between subprocesses that amplify or attenuate signals propagating through component interactions. The standard approach to identify the effects of feedbacks in system responses to forcing is to compare forcing responses in the actual system with responses in a `counterfactual' system in which the feedback process is removed \cite{roe09a}. 
This approach requires high-fidelity simulations of the system behaviors, and even then it is not always clear how to properly ``turn off'' or remove the feedback process because of the interconnectedness of the various sub-processes \cite{merl14a}. 

In addition, the lack of linear superposition makes it very difficult for statistical methods to disentangle the effects of individual mechanisms that interact and intertwine. Together with strong coupling and feedbacks, we require global holistic approaches to parse the complex behaviors generated by nonlinear dynamical systems. For example, how can the clouds of Jupiter simultaneously support apparently random chaotic motion together with coherent structures, like the Great Red Spot, that spontaneously emerge as islands of order in the sea of turbulence? We know the fundamental physics sufficiently well that we can produce such behaviors in simulation, but it is not at all clear how this happens from inspecting the governing equations.

In response to these difficulties, researchers have developed bespoke causal inference methods specifically for nonlinear dynamical systems. The existing toolkit falls largely into two categories, with methods based on either information flows or delay-coordinate embeddings. Information flows \cite{liang05a,ay08,james16a} generalize the early approaches of Wiener \cite{wien56a} and Granger \cite{gran69}, who argued that $Y(t)$ has a causal effect on $X(t)$ if past observations of $Y(t)$ help to better predict the future evolution of $X(t)$. The other class of methods \cite{sugi12a,Harn17a,tsoni18a} utilize the insight of attractor reconstruction using past observations \cite{Pack80,Take81}, called delay-coordinate embeddings in this context. The idea being that if $X(t)$ and $Y(t)$ belong to the same chaotic dynamical system, then their delay embeddings will reconstruct the same geometric object---namely the attractor of their dynamical system. Transfer entropy \cite{schrei00a} and convergent cross mapping \cite{tsoni18a} are the main representatives of information flows and delay embedding methods, respectively. Both have their strengths and limitations, as outlined in e.g. Ref. \cite{camp23a}. 

Here, we introduce a rigorous first-principles approach that is complementary to information flows and delay embedding methods.
Our theory is based on Koopman operators \cite{brun22a} that provide a \emph{global linearization} of general nonlinear dynamics. Rather than track the nonlinear evolution of individual system states, Koopman theory shifts perspectives to the linear evolution of system \emph{observables}---scalar functions of the system state. The global linearization of Koopman operators provides the ability to disentangle the complex web of interactions in general nonlinear systems without restrictive assumptions like small forcing or remaining near equilibrium. 

The power of Koopman operators and their global linearization enables the following advancements for causal discovery in nonlinear dynamical systems. i) We provide a rigorous definition of causal mechanism for dynamical systems, adapting the definition given by the structural causal model framework \cite{bare22a}. ii) This definition highlights the importance of the \emph{time scale} of causal relations in dynamical systems, which is often neglected in other approaches to causality in dynamical systems. iii) We prove an equivalent definition of causal mechanism in terms of Koopman operators that lends itself to a simple data-driven algorithm to quantify causal relations. Whereas most other algorithms for dynamical systems identify causal relations between two univariate time series, our method can naturally accommodate multivariate relations in high-dimensional systems.

Our Koopman theory of causality unfolds as follows. In Section \ref{sec:prelims} we introduce the necessary background and notation for dynamical system flow maps and Koopman operators. Here we define dynamical system components, the analog of variables like $X$ and $Y$ above, over which dynamical causal relations are defined. 
Section~\ref{sec:causality} introduces our definition of causality and causal mechanisms for dynamical systems and its formulation in the Koopman framework. We prove that the two definitions are equivalent. In Section~\ref{sec:dmd_causality} we introduce our data-driven measures to quantify the degree of causal influence and demonstrate them using the coupled \rossler system. Finally, Section~\ref{sec:examples} provides a detailed case study of the Lorenz 96 system and its causal flows driven by advection. 
We show that our data-driven Koopman measure identifies a spatially-asymmetric causal wave that agrees with the wave of information flow identified using a perturbation propagation experiment.
We close with a discussion in Section~\ref{sec:conclusion}.

\section{Preliminaries}
\label{sec:prelims}
\subsection{Dynamical Systems}

Dynamical systems theory describes how a system evolves over time. Mathematically, the system state is given as a point in a \emph{phase space} $\statespace$. The phase space often comes equipped with additional topological or measure theoretic structure depending on the system described. For all the examples we consider in this work, $\M$ is the standard Euclidean space $\Reals^N$, though the results hold equally well for $N$-manifolds. 
Because we will deal almost exclusively with vector quantities, we will not use any notation (e.g. $\vec{\omega}$) to distinguish vector quantities from scalar quantities. Variables should be considered as vectors unless otherwise stated. 

Time-evolution is given by a semigroup of \emph{flow maps}---functions that map the phase space back to itself $\{\Phi^t : \statespace \rightarrow \statespace\}_{t \in \mathcal{T}}$. 
Flow map semigroups are applicable to both discrete time dynamics, with $\mathcal{T} = \mathbb{N}$, and continuous time, $\mathcal{T}= \mathbb{R}_{\geq 0}$. All of our examples are in continuous time, but our formalism is given generally in terms of flow maps that can be in discrete or continuous time. For continuous time, we will 
assume additional structure on the flow maps, namely that they are differentiable (and thus also continuous)~\cite{laso94a}. 

For discrete-time dynamical systems, such as iterated maps~\cite{Miln77,Coll80a}, a generator function $\Phi: \statespace \rightarrow \statespace$ gives the unit-step evolution as 
\begin{align*}
    \state(n+1) = \Phi(\state(n))
    ~.
\end{align*} 
The semigroup elements $\{\Phi^m\}_{m \in \mathbb{N}}$ are then given through $m$-fold composition, so that, e.g. for $m=2$, $\state(n+2) = \Phi^2(\state(n)) = [\Phi \circ \Phi](\state(n))$. 

For continuous time, we consider differential dynamical systems in which $\Phi$ specifies the time-derivative of the system state $\state \in \statespace$
\begin{align*}
    \dot{\state} = \frac{d}{dt} \state = \Phi(\state)
    ~.
\end{align*}
The flow maps $\Phi^t$ are generated from $\Phi$ through time integration, so that 
\begin{align*}
    \state(t) = \Phi^t(\state) = \state + \int_0^t \Phi\bigl(\state(\tau)\bigr) d\tau
    ~.
\end{align*}

For both discrete and continuous time, the time-evolution of the system state $\{\state(t) = \Phi^t(\state)\}_{t \in \mathcal{T}}$ is called an \emph{orbit}, with initial condition $\state$, or simply the orbit of $\state$. 
Note that orbits are functions of both time and initial condition, the latter of which is implicit in our notation. That is, $\omega(t) = \Phi^t(\omega)$ is interpreted as the state at time $t$ starting from initial state $\omega$. To be explicit, this is equivalent to designating the initial state as $\omega_0$ and writing $\omega(t, \omega_0) = \Phi^t(\omega_0)$. We use the simplified notation to avoid excessive clutter. 

As in discrete time, the flow maps form a semigroup through composition so that
\begin{align*}
    \state(t+t') &= \Phi^{t+t'}(\state) = [\Phi^t \circ \Phi^{t'}] (\state)\\ &= \Phi^t \bigl( \Phi^{t'}(\state)\bigr) = \Phi^t(\state(t'))
    ~.
\end{align*}
The infinitesimal generator of the semigroup of flow maps is given by $\Phi$, the time-derivative of the system states, 
\begin{align*}
    \Phi(\state) = \lim\limits_{\tau \rightarrow 0^+} \frac{1}{\tau} \bigl(\Phi^{t+\tau}(\state) - \Phi^t(\state)\bigr) = \frac{d}{dt} \Phi^t(\state) |_{t=0}
    ~.
\end{align*}  

For simplicity, we restrict ourselves to \emph{autonomous} dynamical systems for which $\Phi$ does not have explicit time dependence. 

\subsection{System Components}

Our causal theory concerns how one `part' of a dynamical system may or may not influence the dynamics of another `part'. We now formalize the idea of `parts' of a dynamical system. 
Recall that we consider state spaces $\M$ such that system states $\omega \in \M$ are $N$-dimensional vectors. 
We refer to each scalar element of $\omega$, and equivalently each dimension of $\Omega$, as a \textbf{degree of freedom}. Thus, for the system $\Omega = \mathbb{R}^N$, there are $N$ degrees of freedom. 
Our interest is in causal relations between distinct subsets of degrees of freedom, which we call \emph{components}.

\begin{definition}
    Consider a dynamical system with an $N$-dimensional phase space and hence $N$ degrees of freedom. The degrees of freedom may be partitioned into system \textbf{components}---non-overlapping subsets of degrees of freedom. We denote the partitioning of phase space into component subspaces as $\M = \M_1 \times \M_2 \times \dots \times \M_L$, with $\omega_i \in \M_i$, and write $\omega = [\omega_1 \; \omega_2\; \dots \omega_L]^\top$.
    \label{def:component}
\end{definition}

Note that the partitioning of degrees of freedom into components is mathematically arbitrary and depends on the context of the problem and causality questions being asked. For example, if $\omega = [x \; y \; z]^\top \in \Reals^3$, we can partition into $L = 2$ components with, e.g., $\omega_1 = [x \; y]^\top$ and $\omega_2 = [z]$. We can also choose, e.g., $\omega_1 = [y]$ and $\omega_2 = [x \; z]^\top$ as a valid $L=2$ component partition.

\begin{example}
    To be concrete, we now introduce the model system that will guide our exposition, two \textbf{coupled \rossler oscillators}:
    \begin{align}
        \label{eq:rossler}
        \dot{x}_1 &= -\varphi_1 y_1 - z_1 \nonumber \\
        \dot{y}_1 &= \varphi_1 x_1 + a y_1 + c_1(y_2 - y_1)\nonumber \\
        \dot{z}_1 &= b + z_1(x_1 - d) \nonumber \\
        \; \\
        \dot{x}_2 &= -\varphi_2 y_2 - z_2 \nonumber \\
        \dot{y}_2 &= \varphi_2 x_2 + a y_2 + c_2(y_1 - y_2)\nonumber \\
        \dot{z}_2 &= b + z_2(x_2 - d) \nonumber
        ~,
    \end{align}
    with constants $a,b,c_i,\varphi_i \in \mathbb{R}$.  This is a six-dimensional system, $\Omega = \mathbb{R}^6$, with a natural partitioning into two components $\omega = [\omega_1 \; \omega_2]^\top$ where each component has three degrees of freedom $\omega_j = [x_j \; y_j \; z_j]^\top$. Note that the coupling between the systems in this case is linear, through the $y_i$ degrees of freedom only, with coupling constants $c_i$. Also, \rossler oscillators are weakly nonlinear, with a single quadratic nonlinearity in the $z$ variable dynamics.   
    \label{ex:rossler}
\end{example}

\subsection{Reproducing Kernel Hilbert Spaces}

The Koopman operator, defined below, is a type of map between function spaces.  The most natural type of function space for data analysis are the reproducing kernel Hilbert spaces.  We will need to leverage properties about reproducing kernel Hilbert spaces throughout this work---in particular product spaces of reproducing kernel Hilbert spaces.  

Let $\Omega$ be a non-empty set.  A collection of functions $\F$ mapping $\Omega$ to $\mathbb{C}$ is said to be a \emph{reproducing kernel Hilbert space (RKHS)} if
\begin{itemize}
    \item $\F$ is endowed with the vector space structure over $\mathbb{C}$: $(f+g)(x) = f(x) + g(x)$
    \item $\F$ is equipped with an inner product $\langle \cdot, \cdot \rangle$, with which $\F$ becomes a Hilbert space (complete inner product space)
    \item The evaluation functionals $E_\omega:\F \rightarrow \mathbb{C}$ via $E_\omega(f) = f(\omega)$ are bounded ($\mbox{sup}_{\|f\|\leq 1} |E_\omega(f)| < \infty$) for all $\omega \in \Omega$
\end{itemize}

By the Riesz representation theorem, there exists functions $k_\omega \in \F$ such that $E_\omega(f) = \langle f, k_\omega \rangle$.  These functions are called \emph{the reproducing kernel at $\omega$}. The \emph{kernel function} $\mathsf{K}:\Omega \times \Omega \rightarrow \mathbb{C}$ defined by $\mathsf{K}(\omega, \omega'):= \langle k_{\omega'}, k_{\omega} \rangle = k_{\omega'}(\omega)$ uniquely defines a RKHS by the Moore-Aronszajn theorem~\cite{pauls16a}.  We denote the RKHS associated to the kernel $\mathsf{K}$ by $\mathcal{H}(\mathsf{K})$.  It is known that the reproducing kernels are dense in the Hilbert space:
\[
\mathcal{H}(\mathsf{K}) = \overline{span_{\omega \in \Omega}\{k_\omega\}}
\]
where the overline denotes closure in the norm of $\mathcal{H}(\mathsf{K})$.

Our Koopman causality theory will require reproducing kernels associated to (arbitrary) subdivisions of the phase space $\Omega$ into components. Suppose that there are $L$ components, $\Omega = \Omega_1 \times \Omega_2 \times \dots \times \Omega_L$.  Let $\mathsf{K}_i$ be a kernel function defined over $\omega_i$. Then for $\omega = [\omega_1 \; \omega_2 \; \cdots \; \omega_L]^\top \in \Omega$, one can define the tensor product of functions $[f_1 \otimes f_2 \otimes \dots \otimes f_L](\omega):=f_1(\omega_1)f_2(\omega_1) \dots f_L(\omega_L)$ where $f_i \in \mathcal{H}(\mathsf{K}_i)$.  The tensor product space $\mathcal{H}(\mathsf{K}_1) \otimes \dots \otimes \mathcal{H}(\mathsf{K}_L)$ is a Hilbert space, with inner product 
\[
\langle \otimes f_i, \otimes g_i \rangle = \prod \langle f_i,g_i\rangle_{ \mathcal{H}(\mathsf{K}_i)}
\]
This Hilbert space consists of functions that are limits of sums of these elementary tensors $\otimes f_i$. One may define the \textit{tensor product kernel} $\mathsf{K}: \Omega \times \Omega$ by 
\[
\mathsf{K}(\omega, \omega') = \mathsf{K}_1(\omega_1, \omega_1') \cdots \mathsf{K}_L(\omega_L, \omega_L') 
\]
It is well known that $\mathcal{H}(\mathsf{K})$ is naturally isomorphic to $\mathcal{H}(\mathsf{K}_1) \otimes \dots \otimes \mathcal{H}(\mathsf{\mathsf{K}}_L)$~\cite{pauls16a}. Note that the kernel at a point $\omega'$ for the tensor product kernel is expressed as a tensor product of the individual kernels:
\[
k_{\omega'}(\omega) = \prod_{i=1}^L \mathsf{K}_i(\omega_i,\omega_i') = \prod_{i=1}^L k_{\omega'_i}(\omega_i)
\]
That is, $k_{\omega'} = \otimes_{i=1}^L k_{\omega'_i}$. 
See Appendix~\ref{app:RFFs} for an example using random Fourier features. 

Given that data comes as discrete point-evaluations of functions, RKHS serve as the natural Hilbert spaces for data driven analysis of dynamical systems~\cite{gonz21a}.  While not strictly a requirement, we will simplify the theory discussed in this work by restricting our attention to RKHS's. Any function space $\F$ is thus assumed to be a RKHS unless otherwise specified. 

\subsection{Koopman Operators}

As we will elaborate on further below, nonlinearity in the dynamics $\Phi$ makes it very challenging to track causal dependencies among components of the system state $\state$ over time. Over the last decade, Koopman theory~\cite{brun22a} has emerged as a powerful new approach that provides a \emph{global linearization} of arbitrary dynamical systems, even if $\Phi$ is highly nonlinear. This is in contrast to standard linear stability analysis provided by the Hartman-Grobman theorem~\cite{meiss07a} which states that a dynamical system is locally linear near a hyperbolic fixed point so that the system state evolves linearly in a small neighborhood around the fixed point. 

Of course, if $\Phi$ is nonlinear, the evolution of the system state $\state$ cannot be linear everywhere. Global linearization in Koopman theory is achieved by taking a different perspective. Rather than analyze the evolution of individual system states, Koopman theory analyzes the evolution of system \emph{observables}---scalar functions of the system states $f: \statespace \rightarrow \mathbb{C}$, with $f$ as an element of some function space $\mathcal{F}$. The evolution of observable functions, given by Koopman operators, is linear for all dynamical systems, as long as the function space $\mathcal{F}$ is a vector space. The price that is paid, however, is that the linear evolution of observables generically occurs in infinite dimensions. Throughout, we will assume that our observables arise as functions in an RKHS $\F$ over the phase space $\M$.

The action of a \emph{Koopman operator} $\Koop^t: \mathcal{F} \rightarrow \mathcal{F}$ on an observable $f \in \mathcal{F}$ yields a time-shifted observable, denoted as $f_t \in \mathcal{F}$, that is given by the composition of $f$ with the flow map $\Phi^t$~\cite{Koop31a}
\begin{align}
    f_t(\state) = [\Koop^t f](\state) := [f \circ \Phi^t](\state) = f(\state(t))
    ~.
    \label{eq:koop}
\end{align}
As mentioned, $\Koop^t$ is a linear operator on $\mathcal{F}$, since for all $\alpha \in \mathbb{C}$ and $f,g \in \mathcal{F}$
\[
\Koop^t (\alpha f + g)(\state) = \alpha f(\state(t)) + g(\state(t)) = \alpha \Koop^t f (\state) + \Koop^t g(\state)
\]
Koopman operators inherit the semigroup structure of the flow maps: $\Koop^{t+t'} = \Koop^t  \Koop^{t'}$. The action of the Koopman semigroup $\{\Koop^t\}_{t \in \mathcal{T}}$ on a given observable $f_0$ produces an orbit in function space $\{f_t = \Koop^t f_0\}_{t \in \mathcal{T}}$, analogous to orbits of systems states in the original system. We emphasize that this is an orbit of functions, which can be tricky to work with. However, the evaluation of an observable orbit at a given system state $\state(0)$ yields a scalar time series, given equivalently as
\begin{align}
    \bigl\{f_t\bigl(\state(0)\bigr)\bigr\}_{t\in \mathcal{T}} = \bigl\{f_0\bigl(\state(t)\bigr)\bigr\}_{t \in \mathcal{T}}
    ~.
\end{align}
These observable time series are the objects we work with for data-driven approximation of Koopman operators, discussed in more detail shortly.

Koopman operators over RKHS satisfy the particularly nice property of being bounded.  A linear operator $T: \mathcal{H} \rightarrow \mathcal{K}$ mapping between Hilbert spaces $\mathcal{H}$ and $\mathcal{K}$ is said to be \textit{bounded} if there exists some constant $M$ such that $\|Tf \| \leq M \|f\|$ for all $f \in \mathcal{H}$. It is well known that for linear operators, being bounded is equivalent to being continuous. A Koopman operator whose action on a RKHS stays inside the RKHS is automatically bounded:

\begin{theorem}
    \label{thm:bounded}
    Let $\mathcal{F}$ be a RKHS, and $\Koop^t$ satisfy that $\Koop^t f \in \mathcal{F}$ for all $f \in \mathcal{F}$ ($\Koop^t:\mathcal{F} \rightarrow \mathcal{F}$). Then $\Koop^t$ is bounded (continuous).
\end{theorem}

See Appendix \ref{prf:bounded} for the proof.


\subsection{Finite-Dimensional Approximation: Dynamic Mode Decomposition}

The \emph{dynamic mode decomposition} (DMD) family of algorithms~\cite{will15a, klus16a} provides a matrix approximation $\KoopMat^t$ of the Koopman operator $\Koop^t$ from data using a least-squares regression. The matrix $\KoopMat^t$ converges to a Galerkin projection of $\Koop^t$ on a finite-dimensional subspace of $\F$ in the infinite data limit~\cite{will15a, klus16a}.
Let $\Psi = [\psi_1 \; \psi_2 \; \cdots \; \psi_M]^\top$ be a \emph{dictionary} of scalar-valued observables in $\mathcal{F}$ and $\mathcal{F}_\Psi = span(\Psi) \subset \mathcal{F}$ the finite-dimensional Hilbert subspace spanned by the dictionary $\Psi$.
The \emph{Galerkin projection} of $\Koop^t$ to the finite subspace $span(\Psi) := \F_\Psi \subset \F$ is the unique linear operator $\KoopMat^t: \F_\Psi \rightarrow \F_\Psi$ that satisfies 
\begin{align*}
    \langle \psi_j, \Koop^t \psi_i \rangle = \langle \psi_j, \KoopMat^t \psi_i \rangle, \;\; \text{ for all } \psi_i, \psi_j \in \Psi
    ~.
\end{align*}

We assume that data is taken from a measurable dynamical system $(\M, \Sigma, \mu, \Phi)$ and collected into $D$ pairs $\{\bigl(\omega(n), \omega^t(n)\bigr)\}_D$ of system states $\omega(n)$ and their $t$-shifts $\omega^t(n) := \Phi^t\bigl(\omega(n)\bigr) = \omega(n+t)$. Each $\omega(n)$ is an $N$-dimensional vector. All of our examples below utilize continuous-time dynamical systems with a fixed sampling interval $\delta t$, and the $t$-shifts are given as integer multiples of $\delta t$. 

The data points $\{\bigl(\omega(n), \omega^t(n)\bigr)\}_D$ are collected into data matrices, with the state vectors $\state(n)$ in $\mathbf{\Omega} \in \mathbb{R}^{N \times D}$ and their $t$-shifts $\state^{t}(n)$ in $\mathbf{\Omega}^t \in \mathbb{R}^{N \times D}$:
\begin{align}
    \mathbf{\Omega} = 
    \begin{bmatrix}
        \vline & \vline & \; & \vline \\ 
        \state(n_1) & \state(n_2) & \cdots & \state(n_D) \\
        \vline & \vline & \; & \vline 
    \end{bmatrix}
    \label{eq:data_mat}
\end{align}
and 
\begin{align}
    \mathbf{\Omega}^t = 
    \begin{bmatrix}
        \vline & \vline & \; & \vline \\ 
        \state^{t}(n_1) & \state^{t}(n_2) & \cdots & \state^{t}(n_D) \\
        \vline & \vline & \; & \vline 
    \end{bmatrix}
    ~,
    \label{eq:data_mat_shift}
\end{align}

\noindent respectively.

Using the dictionary functions $\psi_i: \statespace \rightarrow \mathbb{R}$ in $\Psi$, we create observable, or ``boosted'', data matrices 
\begin{align}
    \mathbf{\Psi} = 
    \begin{bmatrix}
        \vline & \; & \vline \\ 
        \Psi(\state(n_1)) & \cdots & \Psi(\state(n_D)) \\ 
        \vline & \; & \vline 
    \end{bmatrix}
    \label{eq:dict_mat}
\end{align}
and
\begin{align}
    \mathbf{\Psi}^{t} = 
    \begin{bmatrix}
        \vline & \; & \vline \\ 
        \Psi(\state^{t}(n_1)) & \cdots & \Psi(\state^{t}(n_D)) \\ 
        \vline & \; & \vline 
    \end{bmatrix}
    ~,
\end{align}
where
\begin{align*}
    \Psi(\state(n)) = [\psi_1(\state(n)) \;  \cdots \; \psi_M(\state(n))]^\top \in \mathbb{R}^M
    ~,
\end{align*} 
and $\mathbf{\Psi}, \mathbf{\Psi}^t \in \mathbb{R}^{M \times D}$. Each column of $\mathbf{\Psi}$ is an $M$-dimensional vector of point-evaluations of each of the $M$ functions $\psi_i$ in the dictionary $\Psi$ at the point $\state(n)$. 

Dynamic mode decomposition seeks the best-fit matrix that evolves the dictionary functions, based on their point-evaluations from the time shifts in the dataset. The DMD matrix $\KoopMat$ minimizes the error
\begin{align}
    &\sum_{n=1}^D || [\Koop^t \Psi](\state(n)) - \KoopMat^t \Psi(\state(n)) ||^2 \nonumber \\
    =&\sum_{n=1}^D || [\Psi \circ \Phi^{t}](\state(n)) - \KoopMat^t \Psi(\state(n)) ||^2 \label{eq:error} \\
    =&\sum_{n=1}^D || \Psi(\state^{t}(n)) - \KoopMat^t \Psi(\state(n)) ||^2 \nonumber
    ~.
\end{align}
The error is equivalent to the Frobenius norm $|| \mathbf{\Psi}^t - \KoopMat^t \mathbf{\Psi}||^2_F$.

The least-squares solution to Eq.~(\ref{eq:error}) is given as 
\begin{align}
    \KoopMat^t := \mathbf{\Psi}^t \mathbf{\Psi}^\dagger
    ~,
    \label{eq:koopmat}
\end{align}
where $\mathbf{\Psi}^\dagger$ is the pseudoinverse of $\mathbf{\Psi}$.

We have emphasized in Eq.~(\ref{eq:error}) that the shifted observable data matrix $\mathbf{\Psi}^t$ is given by point-evaluations of the true action of the full Koopman operator $\Koop^t$. As the Koopman operator is defined in terms of composition with the system dynamics $\Phi$, $\mathbf{\Psi}^t$ is created from the full system dynamics and all its complexity---including relevant causal interactions. 
Said another way, the approximation $\KoopMat^t$ does not involve any simplification of the underlying dynamics $\Phi^t$. 

The DMD approximation $\KoopMat^t$ converges to the Koopman operator $\Koop^t$ in the strong operator topology as $D, M \longrightarrow \infty$~\cite{kord18a}. 

\subsection{Random Fourier Feature Observables}

In this work, we use random Fourier features~\cite{rahi07a} for our dictionary observables as they are simple, flexible, and effective. A \emph{random Fourier feature dictionary} $\Psi_{\text{RFF}}$ is simply a collection of $M$ cosine functions, 

\begin{align}
    \Psi_{\text{RFF}}(\state) = 
    \begin{bmatrix}
        \cos(\phi_1^\top \state + b_1) \\
        \cos(\phi_2^\top \state + b_2)\\
        \vdots \\
        \cos(\phi_M^\top \state + b_M)
    \end{bmatrix}
    ~,
    \label{eq:RFF}
\end{align} 
where each frequency $\phi_j$ is an $N$-dimensional random vector sampled from a specified distribution $p(\phi)$ and $b_j$ are uniformly-distributed random variables.  
Random Fourier features are seen as a Monte Carlo approximation to a positive-definite and shift-invariant kernel, with the distribution $p(\phi)$ determining the kernel~\cite{rahi07a}. For example, random Fourier features with $p(\phi)$ distributed as a spherical Gaussian approximate a Gaussian kernel. We utilize Gaussian $p(\phi)$ with a diagonal covariance matrix for our random Fourier features in the numerical results below. 

In Appendix \ref{app:RFFs}, we recall how random Fourier features approximate shift-invariant kernels and thus the RKHS defined by that kernel. We show that tensor products of random Fourier features are also random Fourier features. This allows us to use RFFs for component function spaces as well as full system function spaces.

\section{Koopman Theory of Causality in Dynamical Systems}
\label{sec:causality}

We now consider causality for nonlinear dynamical systems. First, we define necessary terms and the existence of causal relations in dynamical systems. We then briefly examine the difficulties of quantifying the degree of causal effects in nonlinear systems. This motivates us to formulate causality in the Koopman perspective of the evolution of system observables. The global linearization of observable dynamics makes quantifying causal effects much more tractable for nonlinear systems. We  provide the theory of causality using Koopman operators, then give an operational algorithm based on DMD to quantify causal relations from data. 

\subsection{Causal Mechanism}

To build our theory on solid foundations, we first must address the question: what is a \emph{causal mechanism}? Theories of causal inference generally identify a causal mechanism through a \emph{structural causal model} \cite{pearl09a,bare22a}. A random variable $X$ is assigned values according to its structural assignment: 
\begin{align}
X := m\bigl(\text{pa}(X), \epsilon_X\bigr)
~,
\label{eq:SCM}
\end{align}
where $\text{pa}(X)$ are the causal ``parents'' of $X$ and $\epsilon_X$ is an IID random variable. The map $m$ is interpreted as the causal mechanism of $X$, whose specific value is determined by some set of other variables, $\text{pa}(X)$, that are said to causally influence $X$. The noise term $\epsilon_X$ summarizes all other factors not explicitly modeled, and all noise terms for the full system are assumed to be jointly independent. 
Conceptually, a variable $Y$ is said to be a cause of a variable $X$ if $X$ can change in response to changes in $Y$. 

Notice that time does not appear in the structural causal model in Eqn.~(\ref{eq:SCM}). Typical causal questions addressed in the statistical framework, such as ``Did the aspirin cause my headache to go away?'' or ``Does obesity shorten life?'' \cite{hernan08a}, are not framed in terms of temporal evolution. To include sequential changes over time, statistical methods simply add time dependence to the variables but assume a time-invariance of the qualitative causal relationships \cite{camp23a}. This \emph{causal stationarity} assumption means the functional form of causal mechanisms do not change over time, and thus neither do the causal parents. 

As we will see, causal stationarity is not a reasonable assumption for nonlinear dynamical systems, particularly those with spatial extent. For spatially-extended systems with local interactions, causal influence spreads through the interactions at a finite rate. Therefore, \emph{causal mechanisms can qualitatively change over different time scales}.

\subsection{Causality in Dynamical Systems}



We now adapt the notion of causal mechanism to dynamical systems using flow maps. Because the flow maps are fully parameterized by time, our definition reveals the importance of time scales for causal mechanisms and how they may qualitatively change over time. In addition, our definition highlights the difficulty of identifying causal relations in nonlinear dynamical systems---individual flow maps are complicated objects that generally do not have an interpretable closed-form analytical solution. 

To properly define causal mechanisms for dynamical systems, we first require the following definition:

\begin{definition}
    Let $\M = \M_1 \times \cdots \times \M_L$ be a partition of a dynamical system $(\M, \Phi)$.
    Given an initial state $\omega = [\omega_1 \; \ldots \; \omega_L] \in \M$, recall that the \textbf{orbit} of $\omega$ is $\{\omega(t) := \Phi^t(\omega)\}_{t \in \mathcal{T}}$. 
    We define the \textbf{component orbit} $\omega_j(t)$ as 
    \begin{align}
        \omega_j(t) := P_j \omega(t) = P_j \Phi^t(\omega) \in \Omega_j
        ~,
    \end{align}
    where $P_j$ is the coordinate projection of $\omega$ onto component $\omega_j \in \M_j$.
\end{definition}

The central question for the \emph{existence} of a causal relationship in dynamics is whether the orbit of an ``effect'' component depends on the ``cause'' component. At first this seems limited to a bipartite analysis: e.g., the effect of component $\M_j$ on component $\M_i$. We may, however, want to know the synergistic causal effects that, say, components $\M_j$ and $\M_k$ together have on component $\M_i$. Fortunately, our definition of components is sufficiently flexible to allow for such a multipartite analysis---we can simply combine $\M_j$ and $\M_k$ into a single component $\M_j \times \M_k$. 

For the general exposition of our theory, we introduce the following definitions.

\begin{definition}
    Let $(\M, \Phi)$ be a dynamical system. For causal analysis, we specify a fixed partition
    \begin{align*}
        \M = \M_C \times \M_E \times \M_R
    \end{align*}
    of $\M$ into a \textbf{cause component} $\M_C$, an \textbf{effect component} $\M_E$, and \textbf{remainder component} $\M_R$. The remainder component may or may not be empty. 
\end{definition}


\begin{example}
    Consider the coupled \rossler oscillators in Example~\ref{ex:rossler}.
    We may be interested in the causal influence of component $\omega_1 = [x_1 \; y_1 \; z_1]^\top$ on component $\omega_2 = [x_2 \; y_2 \; z_2]^\top$. In this case, we designate $\M_C = \M_1$, $\M_E = \M_2$, and the remainder $\M_R$ is empty. 
    We may also be interested in more fine-grained details on the causal interaction. For example, the coupling between the two oscillators is through the $y$ degrees of freedom, so we may want to know how much of the causal influence from $\M_1$ to $\M_2$ is carried by $y_1$ versus by $x_1$ and $z_1$. In both of these cases, we let $\M_E = \M_2$. In the first case we would have $\M_C = \{y_1\}$ with $\M_R = \{x_1, z_1\}$, and in the second case $\M_C = \{x_1, z_1\}$ with $\M_R = \{y_1\}$.  
\end{example}

The example above emphasizes that the designation of cause and effect components $\M_C$ and $\M_E$ is arbitrary and depends on the causal questions being investigated. The subscripts $C$ and $E$ are replaced by numerical subscripts, as above, when specifying pre-defined components in a particular system under investigation.

With the cause, effect, and remainder components defined, we can now make precise the existence of causal relations in dynamical systems. 

 \begin{definition}
    \label{def:dyn_cause}
     Let $\M = \M_E \times \M_C \times \M_R$ be defined as above. Fix a $t > 0$ and consider the component flow $P_E \Phi^t$ that is the mapping from a full system state initial condition to the `effect' component $\M_E$ at time $t$
     \begin{align*}
         [\omega_E \; \omega_C \; \omega_R]^\top \mapsto \omega_E(t)
         ~.
     \end{align*}
     We say that $\M_C$ \textbf{dynamically causally influences} $\M_E$ \textbf{at time} $t$ if $\omega_E(t)$ depends on the initial $\omega_C$. We denote this by $\M_C \causes^t \M_E$. Otherwise, $\M_C$ \textbf{does not dynamically causally influence} $\M_E$ \textbf{at time} $t$, and we write $\M_C \ncauses^t \M_E$. 

     If $\M_C \causes^t \M_E$ for any $t \in \mathcal{T}$, we say that $\M_C$ \textbf{dynamically causes} component $\M_E$, and write $\M_C \causes \M_E$. If for all $t$ there is no causal influence, we say that $\M_C$ \textbf{does not dynamically cause} $\M_E$ and write $\M_C \ncauses \M_E$.
 \end{definition}

 \begin{remark}
     The component flow maps $P_E \Phi^t$ play the role of \textbf{causal mechanisms} in our dynamical systems theory, analogous to structural causal models, e.g., Eqn.~(\ref{eq:SCM}). Like other causal inference frameworks, we say that component $\M_C$ has a causal effect on $\M_E$ if and only if $\M_C$ appears in the causal mechanism for $\M_E$. As emphasized in the Introduction, however, the component flow maps are parameterized by the time scale $t$ and, as we will see below, can qualitatively change over time. 
 \end{remark}

\begin{remark}
    Consider two components with $\M = \M_1 \times \M_2 \times \M_R$. If we have that either $\Omega_2 \causes \Omega_1$ or $\Omega_1 \causes \Omega_2$, then $\Omega_1$ and $\Omega_2$ are necessarily components of the same dynamical system. If, on the other hand, we have that both $\Omega_1 \ncauses \Omega_2$ and $\Omega_2 \ncauses \Omega_1$, then $\Omega_1$ and $\Omega_2$ are \textbf{independent} dynamical systems. Finally, if both $\Omega_1 \causes \Omega_2$ and $\Omega_2 \causes \Omega_1$, we say there is a \textbf{dynamical feedback} between the components $\Omega_1$ and $\Omega_2$. 
\end{remark}

\begin{figure*}
    \begin{center}
        \includegraphics[width=1.0 \textwidth]{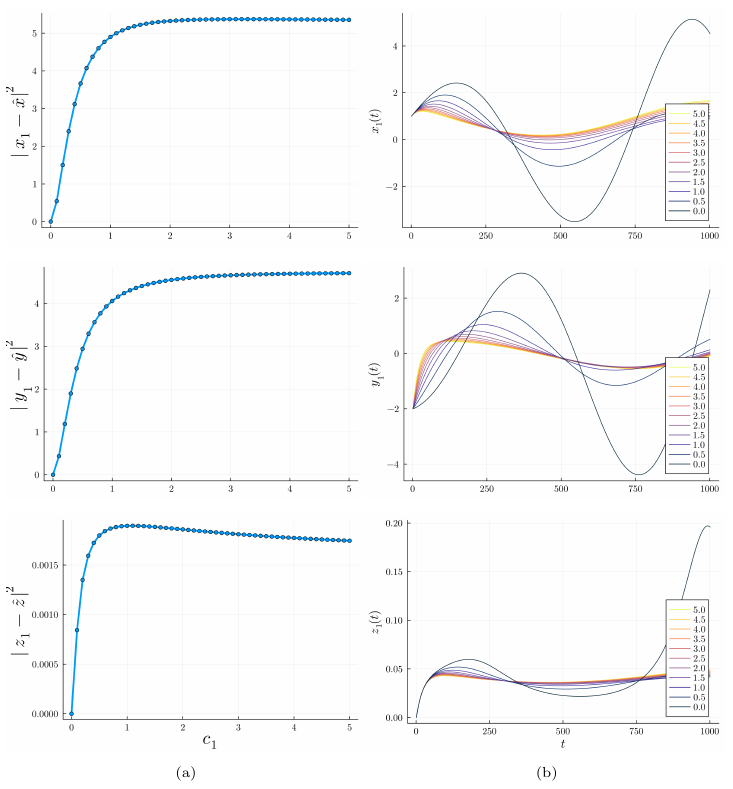}
    \end{center}
    \caption{Some causal phenomenology of the coupled \rossler oscillator system using a counterfactual causality measure. Column (a) shows causality from an asymmetrical system with $\M_2 \causes \M_1$ and $\M_1 \ncauses \M_2$. Column (b) shows sample time series of each degree of freedom with increasing coupling constants $c_1$.}
    \label{fig:rossler}
\end{figure*}

While straightforward to define causality between components in dynamical systems, it is not as obvious how to quantify the \emph{degree of causality} for nonlinear systems. The next example illustrates the issue. 

\begin{example}
    Consider the coupled \rossler system given in Eqn.~(\ref{eq:rossler}) with the two components $\M = \M_1 \times \M_2$. If $c_1 \neq 0$, then $\Omega_2 \causes \Omega_1$ because the evolution of $y_1$ depends on $y_2$. Similarly, if $c_2 \neq 0$ then $\Omega_1 \causes \Omega_2$. If both $c_1 = c_2 = 0$, then the two \rossler oscillators $[x_1 \; y_1 \; z_1]^\top \in \Omega_1$ and $[x_2 \; y_2 \; z_2]^\top \in \Omega_2$ are independent. If $c_1$ and $c_2$ are both nonzero, then there is a dynamical feedback between the two \rossler components of the joint coupled oscillator system. 

    Consider the case of $c_1 \neq 0$ and $c_2 = 0$. It is clear from inspecting the equations of motion that $\M_2 \causes \M_1$ and $\M_1 \ncauses \M_2$. For a specified set of parameters (with $c_1 \neq 0$ and $c_2 = 0$) in the coupled \rossler system, can we identify the degree to which $\M_2 \causes \M_1$? Given that \rossler oscillators are only (weakly) nonlinear in the $z$ degree of freedom and the two \rossler systems are linearly coupled, we might expect the causal influence of $\M_2$ on $\M_1$ to scale with the coupling coefficient $c_1$. 

    To examine causal phenomenology of the coupled \rossler system and explore the interplay of nonlinearity in the oscillators with linear coupling between them, here we use the lens of counterfactuals. Because we are numerically simulating the system and have full control, we can simulate a counterfactual system that is an independent \rossler oscillator with identical parameters and initial conditions to $\M_1$ in the coupled system. The counterfactual system gives what the evolution of $\M_1$ would have been if there were no causal influence from $\M_2$. We can thus define a \emph{counterfactual measure of causality} by computing the mean-squared error between the evolution of $\M_1$ and the counterfactual system: 
    
    \begin{equation}
        \frac{1}{T} \sum_T ||\omega_1(t) - \hat{\omega}_1(t)||^2
        ~,
        \label{eq:counterfactual_measure}
    \end{equation}
    
    \noindent where $\hat{\omega}_1$ is the counterfactual (i.e., a copy of the original system but with $c_1 = 0$). Here, we compare the orbits out to time $T$ that is just below one Lyapunov time. 
    
    We note that this counterfactual measure is really only possible for systems which can be effectively simulated. It is not necessarily possible to compute in \emph{every} simulation. It is simple here because the existence of the causal influence $\M_2 \causes \M_1$ is entirely controlled by whether $c_1$ is zero or not. In more complicated simulations, like climate models, with many nonlinearly-coupled components, it is not straightforward to ``turn off'' causal influence as necessary to create the counterfactual model~\cite{merl14a}. 

    Figure~\ref{fig:rossler} (a) shows the counterfactual causal measure of the $\M_2 \causes \M_1$ coupled \rossler system as a function of the coupling coefficient $c_1$. Each degree of freedom $\{x_1, y_1, z_1\} \in  \M_1$ is analyzed separately in each row. In all cases, we see that causal measure increases sharply as the coupling is first turned on, but quickly saturates to an asymptotic value. Therefore, it is only in the ``weakly coupled'' regime with small $c_1$ that causality is proportional to the coupling constant. As coupling is increased further, it no longer contributes to differences in the orbits between the coupled system and the counterfactual independent system. 

    One may intuit that stronger coupling pushes the system to behave more erratically. Hence the saturation seen in Figure~\ref{fig:rossler} (a) perhaps results from the orbits being pushed to the bounds of the attractor, which limits how far apart an orbit and its counterfactual can be. It turns out that the opposite is true: increased coupling to $\M_2$ acts to stabilize the trajectories of $\M_1$. This is shown in Figure~\ref{fig:rossler} (b), with sample orbits of each degree of freedom given for increasing values of the coupling constant $c_1$. As we see, the orbits with the largest variance are those with small coupling. The causality measure saturates because the orbits converge with increased coupling. Recall that $c_1 = 0$ is the counterfactual case of no coupling. 

    This simple example, with weakly non-linear coupling completely controlled by a single parameter, illustrates the challenges of quantifying causality. While causal influence is, to some degree, controlled by the coupling constants in the coupled \rossler system, it cannot be directly quantified by the coupling constants alone. The effects of the linear coupling also interact with the nonlinearity of the \rossler oscillators. The situation is especially complicated in the presence of feedbacks when both oscillators are coupled together, which we discuss in Appendix \ref{app:rossler_feedback}.
    
    \label{ex:rossler_causal_phenom}
\end{example}

\subsection{Causality in the Koopman Framework}

For simple systems like the coupled \rossler oscillators, the existence of dynamical causal relations is obvious from inspection of the equations of motion---given that they are known. However, for complex systems like a coupled Earth system model, it is not necessarily so clear. More significantly, even for a simple system like the coupled \rossler oscillators with linear coupling and weak nonlinearity, quantifying the degree of causal influence is very challenging, as we have just demonstrated. 

This motivates us to formulate dynamical causality in the Koopman framework. The global linearization in the dynamics of system and component observables turns out to be hugely beneficial for identifying and quantifying causal relations in nonlinear systems directly from data. First, we provide the equivalent of dynamical causality in terms of Koopman operators and flows in function space. The theory then suggests a data-driven algorithm based on DMD to quantify the degree of causal influence with a clear interpretation.

\subsubsection{Component observables and function space embeddings}

Dynamical causality in the Koopman framework manifests as flows in function space. A key step then is to properly formulate the space of observables for individual components as closed subspaces of the function space for the joint dynamical system. 

Let $\Omega = \Omega_1 \times \dots \times \Omega_L$ be a partition of $\Omega$. For each component, let $\mathsf{K}_i$ be a kernel function of $\Omega_i$.  We make the following assumption:

\textbf{Assumption:}  We assume that for each $i$, the RKHS $\mathcal{H}(\mathsf{K}_i)$ contains the constant functions. 

We now define the kernel $\mathsf{K}$ to be the tensor product kernel of $\mathsf{K}_i$, $i=1,\dots,L$, and set $\F = \mathcal{H}(\mathsf{K})$.  Note that by assumption, $\F$ will contain functions which are constant on any arbitrary collection of components of $\Omega$.  This facilitates the following definition:

\begin{definition}
    Let $\Omega= \Omega_1 \times \dots \times \Omega_L$, $\mathsf{K}_i$, $\mathsf{K}$ and $\F$ be as above. Given $X \subset \{1,2, \dots, L\}$, let $\Omega_X = \times_{j \in X} \Omega_j$.  We define $\F_X$ to be the functions in $\F$ which depend only on $\Omega_X$. That is, $\F_X$ is the set of $f \in \F$ such that $f$ restricted to $\Omega_{X^c}$ are constant. 
\end{definition}

\begin{lemma}
    \label{lemma:isometry}
    Let $\M = \M_1 \times \cdots \times \M_L$, $\mathsf{K}_i$, $\mathsf{K}$ and $\F$ be as above.  For each  $X \subset \{1,2, \dots, L\}$,  we have that  $\F_X$  is a  closed subspace in $\F$.
\end{lemma}

See Appendix \ref{prf:isometry} for the proof.

Each of the component subspaces $\F_X$ contains all the information about observables on the space $\M_X$. The function space $\F$ of the full system, on the other hand, contains all the observables on the full state space $\M = \M_1 \times \cdots \times \M_L$. As $\F_X$ are closed subspaces of $\F$, they themselves are (reproducing kernel) Hilbert spaces. We have the following connection between the spaces $\F_X$ and the kernel function $\mathsf{K}$:

\begin{proposition}
\label{prop:dense span}
    Let $X \subset \{1,2, \dots, L\}$, and $\mathsf{K}_i$, $\F_X$ be as above.  Define the kernel function $\mathsf{K}_X: \Omega \times \Omega \rightarrow \mathbb{C}$ by 
    \[
     \mathsf{K}_X(\omega', \omega) := \prod_{j \in X} \mathsf{K}_j(\omega_j', \omega_j)
    \]
    for all $\omega, \omega' \in \Omega$.  Define the functions $k_{X,\omega}: \Omega \rightarrow \mathbb{C}$ by $k_{X,\omega}(\omega') = \mathsf{K}_X(\omega', \omega)$. Then
    \[
    \F_X  = \overline{span_{\omega \in \Omega }} \{k_{X,\omega} \}
    \]
\end{proposition}

See Appendix \ref{prf:dense span} for the proof.

Therefore, the component kernel functions $\mathsf{K}_X$ define the component function subspaces $\F_X = \mathcal{H}(\mathsf{K}_X)$ through their reproducing kernels $k_{X,\omega'}$. The component reproducing kernels and function subspaces naturally compose to the full joint system through tensor products.

\subsubsection{Existence of causal relations from Koopman operators}

With the component observable function subspaces appropriately defined, we now turn to identifying the Koopman equivalent of the existence of dynamical causality. The central idea of our theory is that the existence and absence of dynamical causal relations produce distinct behaviors in the flow of functions in the observable space $\F$. 

Recall from Eqn.~(\ref{eq:koop}) above that the action of a Koopman operator on an observable $f \in \F$ gives a new function in $f_t \in \F$ defined by the composition of $f$ with the flow map $\Phi^t$. Therefore, the time-shifted function $f_t$ inherits functional dependence from the flow map. Intuitively, 
if component $\M_C$ has no causal influence on component $\M_E$, then the Koopman evolution of functions $\F_C$ that initially only depend on $\M_C$ should be indifferent to values of $\omega_E \in \M_E$. 
From Definition~\ref{def:dyn_cause}, if $\M_C \ncauses \M_E$, then the component flow map for $\M_E$ does not depend on $\M_C$ and so the action of Koopman operators on functions in $\F_E$ will inherit this functional independence. We now make these ideas precise. 

\begin{definition}
    Let $(\M, \Phi)$ be a dynamical system as above. Consider a fixed partition
    \[
    \M = \M_C \times \M_E \times \M_R
    \]
    of $\M$ into a `cause', `effect', and `residual' components respectively.  As above, define $\F_X$  as the subset of functions in $\F$ that are constant in all components other than $\M_X$, for $X \subset \{C,E,R\}$. We say that $\M_C$ \textbf{does not Koopman causally influence} component $\M_E$ at time $t$, if $\Koop^t \F_E \subset \F_{E,R}$.  We denote this by $\M_C \ncauses^t_\Koop \M_E$. Otherwise, $\M_C$ \textbf{Koopman causally influences} component $\M_E$ at time $t$, and we denote this by $\M_C \causes^t_\Koop \M_E$.  If for all $t$ we have that $\M_C$ does not Koopman causally influence component $\M_E$, we say that $\M_C$ does not Koopman cause $\M_E$, and write $\M_C \ncauses_\Koop \M_E$
    \label{def:koop_cause}
\end{definition}

To make this definition explicit; the condition $\Koop^t \F_E \subset \F_{E,R}$ can be phrased in words as follows.  Let $f$ be a function that is constant on (at least) the component $\M_C$.  If the Koopman-evolved function $\Koop^t f$ still does not depend on the component $\M_C$, i.e. $\Koop^t f \in \F_{E,R}$, then there is no causality from component $\M_C$ to $\M_E$ as far as the Koopman operator can detect at time $t$. Thus, $\M_C \ncauses^t_\Koop \M_E$. As a shorthand, the lack of a $C$ subscript in $\F_{E,R}$ means there is no dependence on component $\M_C$.

This definition is the equivalent of dynamical causality in the Koopman framework. 

\begin{theorem}
    \label{thm:causal_equiv}
    Dynamical causal influence and Koopman causal influence are equivalent, 
    \begin{align*}
        \M_C \causes^t_\Koop \M_E \iff \M_C \causes^t \M_E
    \end{align*}
    Correspondingly, 
    \begin{align*}
        \M_C \ncauses^t_\Koop \M_E \iff \M_C \ncauses^t \M_E
    \end{align*}
\end{theorem}

See Appendix \ref{prf:causal_equiv} for the proof.

\section{Quantifying Causal Influence From Data with DMD Koopman Approximations}
\label{sec:dmd_causality}

The equivalence of dynamical causality and Koopman causality, given by Theorem~\ref{thm:causal_equiv}, opens new avenues to assess causality in nonlinear systems. The causal flow of observables that identify Koopman causality in Definition~\ref{def:koop_cause} suggests a simple adaptation of the DMD algorithm to \emph{quantify} causal relationships directly from data. 

From Definition \ref{def:koop_cause}, component $\M_C$ has a causal effect on $\M_E$ if the Koopman evolution of functions in the component subspace $\F_E$ induces a dependence on component $\M_C$ in the resulting time-shifted functions. We can thus test how well $\Koop^t \F_E$ can be represented in a basis of functions in the joint observable subspace $\F_{E,C}$. Said another way, we start with functions that depend only on component $\M_E$ and then test how much dependence on component $\M_C$ is added, or not, after evolving with a Koopman operator.

Intuitively, this is a way to quantify the unique causal contribution of $\M_C$ on the evolution of component $\M_E$. If there are additional causal contributions from the remainder component $\M_R$, then $\Koop^t \F_E$ cannot be completely represented in the basis $\F_{E,C}$ since there are missing dependencies from $\M_R$. However, if $\M_C$ encompasses the totality of causal influence on the evolution of $\M_E$, the $\Koop^t \F_E$ can be fully represented in the joint basis $\F_{E,C}$. Note though, that if $\M_E$ is an independent dynamical system with no additional causal influences, then $\Koop^t \F_E$ will also be fully represented by the joint space $\F_{E,C}$ because $\Koop^t \F_E \subseteq \F_E$ and $\F_E \subset \F_{E,C}$. Therefore, we need to add a baseline comparison of how well $\Koop^t \F_E$ can be represented by functions in the marginal state space $\F_E$ that depend only on the `effect' component $\M_E$.

\subsection{Data-Driven Causality Measures}

With these insights, we now define a data-driven causality measure based on the DMD algorithm.  This is done by fitting two different DMD models to evolve test functions in $\F_E$: a \emph{joint model} that uses basis functions in $\F_{E,C}$ and a \emph{marginal model} that uses basis functions in $\F_E$ only. By measuring the difference in these two models, we derive an empirical measure of nonlinear causality through the Koopman operator.

Our Koopman causality algorithm requires data from the chosen `cause' and `effect' components. The data must come as triples $(\omega_E(n), \omega^t_E(n), \omega_C(n))$, where $\omega_E^{t} := P_E \Phi^{t} ([\omega_E \; \omega_C \; \omega_R]^\top)$ is the time-shift of $\omega_E$ and $n$ is the time index of the data. We split the dataset into $\texttt{train}$ and $\texttt{test}$ sets. 

Because we are testing how well functions $\F_E$ can be represented in the joint space $\F_{E,C}$ compared to the marginal space $\F_E$ after Koopman evolution, the algorithm additionally requires separate dictionaries of basis functions for each of these spaces. For the input \emph{test functions} $f_E \in \F_E$, we have found it best to use the identity observables $f_E(\omega_E) = \omega_E$ in practice. Their efficacy is likely because the identity observables are generally not in a finite invariant subspace~\cite{brun16a}, and so adding additional basis functions tends to improve how well the model predicts these observables.



For both the joint and marginal models, we include the test functions (i.e. the identity observables) in the DMD function dictionary $\Psi$, along with the additional basis functions with which we want to express their time-shifts. 
We use random Fourier features for these additional basis functions---see Section~\ref{sec:prelims} E.1. above. 
The marginal dictionary functions are given as
\begin{align}
    \psi_E(\omega_E) = \cos(\phi_E^\top \omega_E + b)
    ~,
    \label{eq:rff_marg}
\end{align}
and the joint dictionary functions are 
\begin{align}
    \psi_{E,C}(\omega_E, \omega_C) = \cos(\phi_{E,C}^\top [\omega_E \; \omega_C]^\top + b)
    ~.
    \label{eq:rff_joint}
\end{align}

Using data matrices as defined above for dynamic mode decomposition, the \textbf{marginal model} is defined as 
\begin{align}
    \KoopMat^t_{\text{marg}} = 
    \mathbf{\Omega}_E^{t} 
    \begin{bmatrix}
        \mathbf{\Omega}_E \\
        \mathbf{\Psi}_E
    \end{bmatrix}^\dagger
    ~,
    \label{eq:marginal_model}
\end{align}
where $\mathbf{\Psi}_E$ is the data matrix formed by applying the marginal dictionary of random Fourier features $\psi_E(\omega_E)$ to the training data matrix $\mathbf{\Omega}_E$.
The \textbf{joint model} is similarly defined as 
\begin{align}
    \KoopMat^t_{\text{joint}} = 
    \mathbf{\Omega}_E^{t} 
    \begin{bmatrix}
        \mathbf{\Omega}_E \\
        \mathbf{\Psi}_{E,C}
    \end{bmatrix}^\dagger
    ~,
    \label{eq:joint_model}
\end{align}
with the data matrix $\mathbf{\Psi}_{E,C}$ formed from the joint dictionary of random Fourier features $\psi_{E,C}(\omega_E, \omega_C)$ applied to the columns of $\mathbf{\Omega}_E$ and $\mathbf{\Omega}_C$. 
The only difference between the two models is the additional dictionary functions of $\mathbf{\Psi}_E$ for the marignal model and $\mathbf{\Psi}_{E,C}$ for the joint model. 

Note that we only care about the approximated Koopman evolution of the test functions in $\F_E$, hence the DMD models $\KoopMat^t_{\text{marg}}$ and $\KoopMat^t_{\text{joint}}$ are not square. For given values of $\omega_E \in \M_E$ and $\omega_C \in \M_C$, the point-evaluation of the time-shifted test function $[\Koop^{t} f_E]([\omega_E \; \omega_C \; \omega_R]^\top) := \omega_E^{t}$ is approximated by the marginal model as 
\begin{align}
    \widetilde{\omega}_E^{t} |_{\text{marg}} := \KoopMat^t_{\text{marg}}
    \begin{bmatrix}
        \omega_E \\
        \Psi_E(\omega_E)
    \end{bmatrix}
    ~,
    \label{eq:marg_pred}
\end{align}
and by the joint model as 
\begin{align}
    \widetilde{\omega}_E^{t} |_{\text{joint}} := \KoopMat^t_{\text{joint}}
    \begin{bmatrix}
        \omega_E \\
        \Psi_{E,C}(\omega_E, \omega_C)
    \end{bmatrix}
    ~.
    \label{eq:joint_pred}
\end{align}

For both cases, model error is evaluated over the set of \texttt{test} points $\{\omega_E, \omega_E^{t}, \omega_C\}_{\texttt{test}}$
\begin{align}
    \text{model error } = \frac{1}{N_{\texttt{test}}}\sum_{\texttt{test}} || \widetilde{\omega}_E^{t} - \omega_E^{t}||^2
    ~,
    \label{eq:model_error}
\end{align}
with the marginal model error given using $\widetilde{\omega}_E^{t} = \widetilde{\omega}_E^{t}|_{\text{marg}}$ and the joint model error given using $\widetilde{\omega}_E^{t} = \widetilde{\omega}_E^{t}|_{\text{joint}}$.

We can now define our \textbf{data-driven measure of causality}, denoted $\M_C \ddcauses \M_E$, as
\begin{align}
    \M_C \ddcauses \M_E := \text{marginal error}-\text{joint error}
    ~.
    \label{eq:ddcause}
\end{align}

\begin{remark}
    Our causality measure is conceptually similar to Wiener-Granger causality~\cite{wien56a,gran69} in that we compare the predictive power of a marginal model to a joint model. In fact, if delay-coordinate embedding observables are used for the additional dictionary functions $\mathbf{\Psi}$ (known as delay DMD or Hankel DMD \cite{arba17a}) the marginal and joint model predictions in Eqns. (\ref{eq:marg_pred}) and (\ref{eq:joint_pred}) are functionally equivalent to those of Wiener-Granger causality---predictions are linear functions of past observables. Note though that the evaluation of model error in Wiener-Granger causality is typically different from Eqn. (\ref{eq:model_error}), typically comparing the variance of model residuals. 
    
    Although conceptually similar, we emphasize that \textbf{our method is not an extension of Wiener-Granger causality}. First, our method is rigorously derived starting from a clear definition of causal mechanism in dynamical systems, Definition~\ref{def:dyn_cause}. Second, while Granger causality assumes a linear relationship between delay embeddings of system components, our method utilizes the linear \emph{flow of observables} that holds generally for nonlinear dynamical systems. Third, and most important, our method explicitly \textbf{does not utilize time delay embeddings}. As detailed in, for example, Reference~\cite{rupe22b}, delay embeddings of one system component encodes predictive information from all other system components. Thus, marginal models in Granger causality already contain predictive information from other components, limiting its ability to distinguish the influence of the added component in the joint model.
    
\end{remark}

\begin{figure}
    \begin{center}
        \includegraphics[width=0.5 \textwidth]{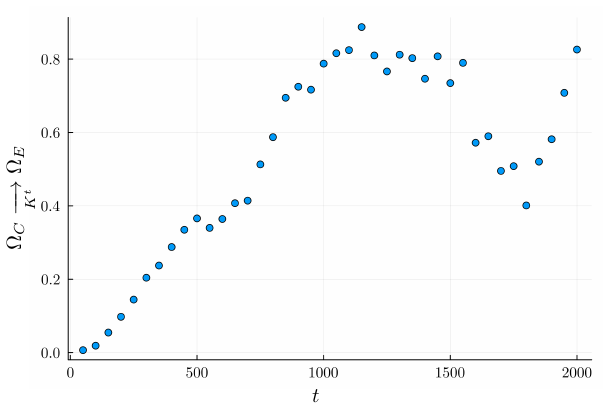}
    \end{center}
    \caption{Causal measure per time of coupled \rossler.}
    \label{fig:rossler_causes_per_time}
\end{figure}

\begin{example}
    In Figure~\ref{fig:rossler_causes_per_time} we show $\M_2 \ddcauses \M_1$ for many time shifts $t$ of the coupled \rossler system with nonlinearity $d=5.7$ and coupling constants $c_1 = 0.5$ and $c_2 = 0.0$ (and so $\M_2 \causes \M_1$ and $\M_1 \ncauses \M_2$). The coupled \rossler system is continuous-time and we use an integration step $\delta t = 0.01$ for data generation. Because the system is continuous-time, a small time shift $t$---given here as the number of integration steps---imparts a small causal dependence of $\M_2$ onto $\M_1$. This is what we see in Figure~\ref{fig:rossler_causes_per_time} with a roughly linear dependence of $\M_2 \ddcauses \M_1$ on $t$ for $t \leq 500$. The causal dependence continues to grow until about $t = 1000$ where it plateaus. Interestingly, around $t = 1500$, the causal dependence actually decreases momentarily, before increasing again back up to the plateau value by $t = 2000$. The Lyapunov time for this system is just under $t = 1000$ time steps. 
    
    We emphasize again that while it is clear from the equations of motion that $\M_2 \causes \M_1$, and it is intuitive that the causal effect should be low at short time scales, it is not at all obvious from inspecting the equations how the causal influence should behave at longer time scales. In Example~\ref{ex:rossler_causal_phenom}, we emphasized the interplay between component coupling and nonlinearity in each components' dynamics and how these combine in non-intuitive ways to produce causal effects. Nonlinear systems, like the \rossler oscillator, can experience intermittent bursting behaviors over longer time scales, which may contribute to the decreased causal influence we see between $t = 1500$ and $t = 2000$. 
    \label{ex:rossler_cause_per_time}
\end{example}

\subsection{Causality Over Time with Conditional Forecasting}

The causal measure in Eqn.~(\ref{eq:ddcause}) is the core data-driven tool to apply our Koopman theory of causality in practice. By definition, it determines causal influence \emph{at the given time scale} $t$. One may be interested more holistically in how the causal influence behaves over time. Because the causal measure in Eqn.~(\ref{eq:ddcause}) has to be computed at each time scale $t$, it can be costly to compute the temporal behavior of causal relations. 

\begin{figure}[h]
    \begin{center}
        \includegraphics[width=0.48 \textwidth, trim={0.0cm 5.2cm 4.6cm 0.0cm}, clip]{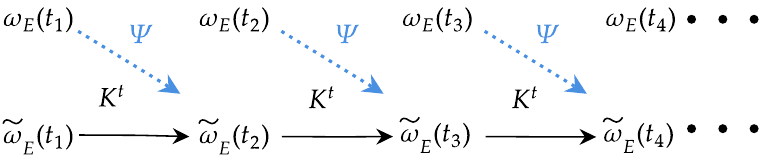}

    \end{center}
    \caption{Graphical depiction of a marginal conditional forecast for `effect' component $\M_E$. Predicted values $\widetilde{\omega}_E(t_n)$ are fed back into the model, shown with solid black arrows, while values $\omega_E(t_n)$ from the \texttt{test} data are plugged into the additional dictionary functions $\Psi$, as shown with dotted blue arrows.}
    \label{fig:causalDMD}
\end{figure}

\begin{figure*}
    \begin{center}
        \includegraphics[width=1.0 \textwidth]{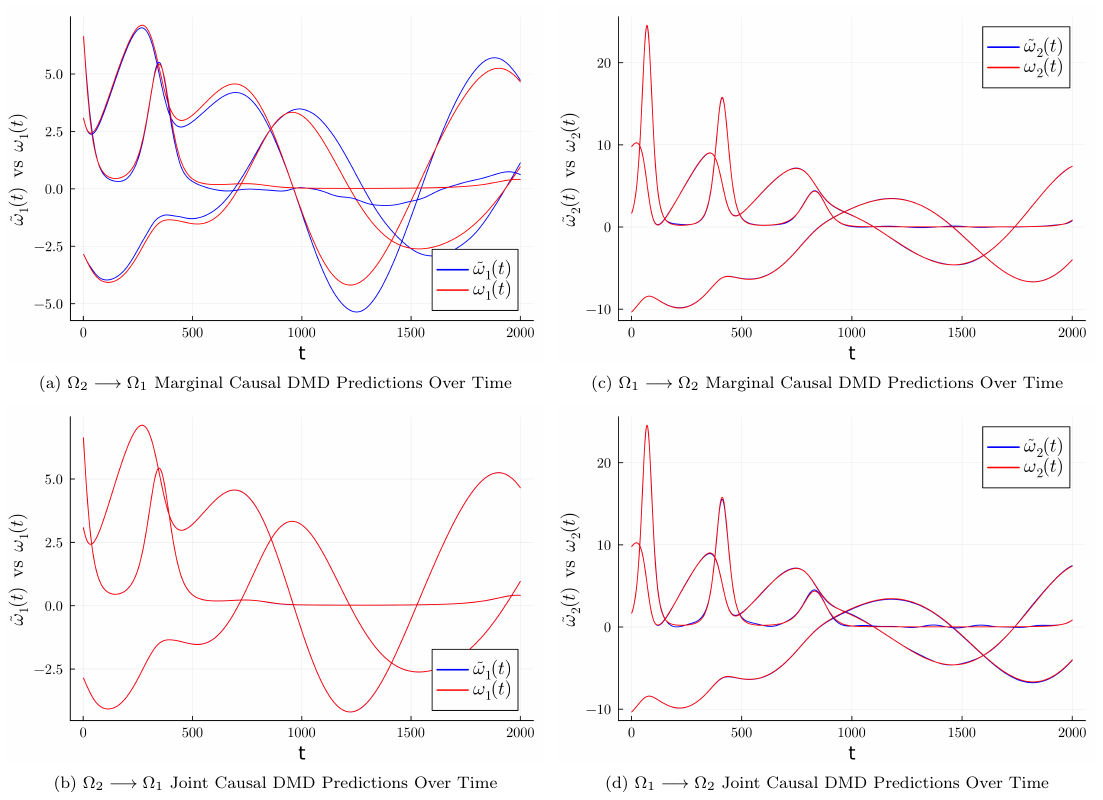}
    \end{center}
    \caption{Conditional forecasts for $\M_2 \causes \M_1$, $\M_1 \ncauses \M_2$ coupled \rossler system.}
    \label{fig:rossler_causalDMD}
\end{figure*}

If the time shift $t$ in the \texttt{training} set is small and the \texttt{test} set is a time series (sampled at the small shift $t$), then we can assess causal relations over time much more efficiently using \emph{conditional forecasting}. The term `forecasting' is used because the algorithm sequentially feeds its output back into its input. However, we emphasize that these models \emph{cannot} predict the evolution of the test functions for times outside the \texttt{test} set---the forecasts are `conditional' because they require inputs from the \texttt{test} set. Of course, the goal here is to assess causality, and not to predict how the observables will evolve in the unseen future---standard DMD and its many variants are used for that. 



\textbf{Conditional forecasting} utilizes the same form of marginal and joint models given in Eqns. (\ref{eq:marginal_model}) and (\ref{eq:joint_model}), respectively, using the small time shift $t$.
The models use dictionaries that consist of two pieces $[\omega_E \; \Psi]^\top$ with $\Psi = \Psi_E(\omega_E)$ for the marginal model and $\Psi = \Psi_{E,C}(\omega_E, \omega_C)$ for the joint model. For both cases, we sequentially forecast the first piece (the identity observable test function), while the second piece, $\Psi$, is always evaluated using the corresponding values of $\omega_E$ and $\omega_C$ from the \texttt{test} set---see Figure~\ref{fig:causalDMD}. For the marginal case, we actually could use the predicted values of $\omega_E$ to evaluate the dictionary functions $\Psi_E(\omega_E)$---this is just standard DMD. However, causal evaluation requires comparison of the marginal and joint models, so we use the same conditional forecasting procedure for both. We emphasize again that this `forecast' only approximates the \texttt{test} time series $\{\omega_E\}_{\texttt{test}}$.

Given initial values $(\omega_E(t_0), \omega_C(t_0))$ from the \texttt{test} data, the conditional forecast proceeds as follows. 
The first prediction is made using the initial \texttt{test} values as 
\begin{align*}
    \widetilde{\omega}_E(t_1)|_{\text{marg}} = \KoopMat^t_{\text{marg}}
    \begin{bmatrix}
        \omega_E(t_0) \\
        \Psi_E(\omega_E(t_0))
    \end{bmatrix}
    ~,
\end{align*}
for the marginal model and 
\begin{align*}
    \widetilde{\omega}_E(t_1)|_{\text{joint}} = \KoopMat^t_{\text{joint}}
    \begin{bmatrix}
        \omega_E(t_0) \\
        \Psi_{E,C}(\omega_E(t_0), \omega_C(t_0))
    \end{bmatrix}
    ~,
\end{align*}
for the joint model.
The next time step is predicted as 
\begin{align*}
    \widetilde{\omega}_E(t_2)|_{\text{marg}} = \KoopMat^t_{\text{marg}}
    \begin{bmatrix}
        \widetilde{\omega}_E(t_1) \\
        \Psi_E(\omega_E(t_1))
    \end{bmatrix}
    ~,
\end{align*}
for the marginal model and 
\begin{align*}
    \widetilde{\omega}_E(t_2)|_{\text{joint}} = \KoopMat^t_{\text{joint}}
    \begin{bmatrix}
        \widetilde{\omega}_E(t_1) \\
        \Psi_{E,C}(\omega_E(t_1), \omega_C(t_1))
    \end{bmatrix}
    ~,
\end{align*}
for the joint model. Note that the previously predicted value $\widetilde{\omega}_E(t_1)$ is fed back in for the identity piece of the dictionary for both models.
This procedure continues in the same way for the duration of the \texttt{test} time series, with predicted values being fed back into the first piece of the model dictionary and the second additional dictionary elements always evaluated with values directly from the test data. 

The marginal and joint model errors can be evaluated, similar to the counterfactual causality measure in Eqn.(\ref{eq:counterfactual_measure}), by taking the mean-square error between the \texttt{test} series and conditional forecasts. One benefit of the conditional forecasting approach is that it provides a visualization of how well the marginal and joint models perform in predicting the \texttt{test} time series, as we now demonstrate with the \rossler system.


\begin{example} 
    Figure~\ref{fig:rossler_causalDMD} shows conditional forecasts for the $\M_2 \causes \M_1$ coupled \rossler system.
    The left column shows forecasts for $\M_2 \causes \M_1$ and the right for $\M_1 \causes \M_2$. 

    In the left column, we see non-zero prediction error from the marginal model, even though values from the \texttt{test} data are used at each step to evaluate $\Psi_1(\omega_1)$. The error occurs because the marginal model has no information from component $\M_2$ that causally influences $\M_1$.
    When this information from $\M_2$ is incorporated in the joint model, it is able to predict $\M_1$ over at least $2000$ time steps (with integration step size $\delta = 0.01$) with very little error. Because the remainder component here is empty and $\M = \M_1 \times \M_2$, the `cause' component $\M_C = \M_2$ captures all causal influences on the `effect' $\M_E = \M_1$.     
    As also observed in Figure~\ref{fig:rossler_causes_per_time}, the marginal model is able to predict reasonably well for several hundred time steps before significant error occurs. For continuous-time systems, causal influence takes some time to propagate. 

    The right column switches roles of the components, with $\M_C = \M_1$ and $\M_E = \M_2$ to test whether component $\M_1$ causally influences $\M_2$, which we know it does not since $c_2 = 0$. As expected then, we see essentially no difference in the predictions of the marginal and joint models. Because there are no external causal influences on $\M_2$, the marginal function space $\F_2$ of component $\M_2$ is Koopman-invariant $\Koop^t \F_2 \subseteq \F_2$. The evolution of test functions are thus well-expressed in the marginal basis $\F_2$ (which also requires that the RFF basis functions cover the space reasonably well).
\end{example}

\subsection{Comparison with Other Causal Measures}
We now compare our Koopman causality meausure with other common causality methods on the uni-directional \rossler model with $\M_2 \causes \M_1$ and $\M_1 \ncauses \M_2$. 
Concretely, we will contrast the results of Koopman causality, Granger causality (GC), and Cross Convergent Mapping (CCM) algorithms. See Appendix \ref{app:comparison} for more details on these methods and the comparison tests. For this comparison, the coupled \rossler system is run 100 times with different initial conditions chosen randomly on the attractor. Causality scores are then measured using Koopman causality, GC and CCM algorithms.  Table \ref{tab:comparison} reports the detection rate of each algorithm: 

\begin{table}[h]
\begin{tabular}{|l|l|l|l|l|}
\hline
\textbf{Method} & $ \Omega_1 \rightarrow \Omega_2$ & $ \Omega_2 \rightarrow \Omega_1$ & Both & Neither \\ \hline
\textit{GC} & 0.00  & 0.00  & 1.00 & 0.00    \\ 
\hline
\textit{CCM} & 0.01 & 0.01 & 0.96 & 0.02     \\ 
\hline
\textit{Koopman} & \textbf{0.00} & \textbf{1.00} & \textbf{0.00} & \textbf{0.00}\\ 
\hline
\end{tabular}
\caption{Comparison of different causality methods at detecting causality in the uni-directional coupled \rossler system with $\M_2 \causes \M_1$. \label{tab:comparison}}
\end{table}

Because of the uni-directional coupling of our test system, the methods should find that $\M_2 \causes \M_1$ and $\M_1 \ncauses \M_2$. As we see in Table \ref{tab:comparison}, GC and CCM identify that there is a causal relation between the components, but they incorrectly identify a bi-directional causal relation. Our Koopman causality measure, on the other hand, identifies the correct uni-directional causal relationship in all cases. 

\section{Confounding}
\label{sec:confounding}
\begin{figure}
    \begin{center}
        \includegraphics[width=0.5 \textwidth, trim={0cm 0cm 0.8cm 0.8cm}, clip]{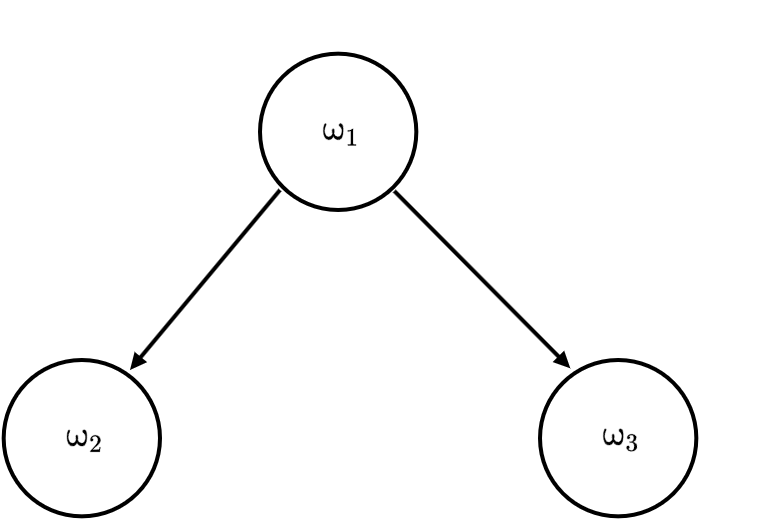}
    \end{center}
    \caption{Causal graph depicting $\omega_1$ as a confounder for components $\omega_2$ and $\omega_3$}
    \label{fig:confounder}
\end{figure}

\emph{Confounding} is a core concept in causal inference since it leads to spurious associations that are not causal. In the dynamical systems literature, confounders are sometimes also called \emph{common drivers}~\cite{rung19a}. The simplest case is shown in Figure~\ref{fig:confounder} with $\omega_1$ acting as a confounder for $\omega_2$ and $\omega_3$. There is no causal influence between $\omega_2$ and $\omega_3$, but they become associated (e.g. correlated) due to their common driver $\omega_1$ that has causal influence on each of them individually. 

For nonlinear dynamical systems, which are typically highly interconnected, confounding is inescapable. This is a primary reason why causality is so difficult for nonlinear systems. For statistical approaches to causal inference that identify causation from association, it is necessary that all confounding variables are observed so that they can be accounted for using the do calculus~\cite{pearl09a}. Although we emphasize again that these approaches cannot handle cyclic causal relations, which are ubiquitous in nonlinear systems. 

The challenges of confounding, particularly for nonlinear dynamical systems, highlights the power of our approach. In particular, it is based directly on the dynamical causal mechanism in Definition~\ref{def:dyn_cause} and does not rely on associations. For a dynamical system corresponding to Figure~\ref{fig:confounder}, components $\omega_2$ and $\omega_3$ are associated due to their confounder $\omega_1$, but they are not causally coupled because $\omega_3$ does not appear in the equations of motion for $\omega_2$ and vice versa. Varying one of them, with $\omega_1$ fixed, will not effect the other. 
Therefore, \textbf{our theoretical framework is indifferent to the presence of confounders, whether they are observed or not}. 

That said, while our theoretical framework naturally handles confounding, our numerical algorithm to approximate causal influence utilizes DMD, which is constructed via regression---a form of association. Thus, it is necessary to test the ability of our numerical algorithm to handle confounding. To this end, we introduce the \textbf{\rossler confounder} model with the causal structure of Figure~\ref{fig:confounder},
\begin{align}
    \label{eq:confounder}
    \dot{x}_1 &= -\varphi_1 y_1 - z_1 \nonumber \\
    \dot{y}_1 &= \varphi_1 x_1 + a y_1 \nonumber \\
    \dot{z}_1 &= b + z_1(x_1 - d) \nonumber \\
    \; \nonumber\\
    \dot{x}_2 &= -\varphi_2 y_2 - z_2 \nonumber \\
    \dot{y}_2 &= \varphi_2 x_2 + a y_2 + c_2(y_1 - y_2)\\
    \dot{z}_2 &= b + z_2(x_2 - d) \nonumber\\
    \; \nonumber\\
    \dot{x}_3 &= -\varphi_3 y_3 - z_3 \nonumber \\
    \dot{y}_3 &= \varphi_3 x_3 + a y_3 + c_3(y_1 - y_3)\nonumber \\
    \dot{z}_3 &= b + z_3(x_3 - d) \nonumber
    ~.
    \end{align}

For this model, we partition it as usual into components $\M_1 = \{x_1, y_1, z_1\}$, $\M_2 = \{x_2, y_2, z_2\}$, and $\M_3 = \{x_3, y_3, z_3\}$. Written more compactly as 
\begin{align*}
    \dot{\omega}_1 &= f_1(\omega_1) \\
    \dot{\omega_2} &= f_2(\omega_2, \omega_1) \\
    \dot{\omega_3} &= f_3(\omega_3, \omega_1)
    ~,
\end{align*}
we see that the \rossler confounder model in (\ref{eq:confounder}) has the causal structure of Figure~\ref{fig:confounder}. Note that the causal relationships between components $\M_1$ and $\M_2$ are the same as the uni-directional coupled \rossler model considered above, since component $\M_3$ has no causal influence over $\M_1$ or $\M_2$. Components $\M_1$ and $\M_3$ similarly behave like a uni-directional coupled \rossler model. Therefore, our primary interest in the \rossler confounder is the relationship between $\M_2$ and $\M_3$. From the equations of motion, we see that neither have causal influence on the other, $\M_2 \ncauses \M_3$ and $\M_3 \ncauses \M_2$---they are causally independent. But, they are not statistically independent; the common driving of the confounder $\M_1$ adds non-causal association between $\M_2$ and $\M_3$. 


To test our method in the presence of confounding, we perform a similar experiment as above using two different configurations of the \rossler confounder model, each with 100 test runs. 

\subsection{Symmetric coupling}
First, we look at the case of symmetric coupling, with $c_2 = c_3$ (note that $\M_2$ and $\M_3$ have different initial conditions, so they are not merely copies of one another). As seen in the \textit{Symmetric} row of Table~\ref{tab:confouder}, our method correctly identifies no causal relationship between $\M_2$ and $\M_3$ in most of the test cases. However, unlike the simpler coupled \rossler case in Table~\ref{tab:comparison}, our method fails in some cases here due to the spurious association induced by the confounder $\M_1$. 

It is remarkable that our method successfully identifies no causal relation in so many cases out-of-the-box, despite the strong statistical association between $\M_2$ and $\M_3$. Further, it does so without any information from the confounder $\M_1$ that induces the statistical association. That is, $\M_1$ is an \emph{unobserved} confounder in this experiment. Although our theoretical framework naturally handles confounding, it is significant to see this demonstrated so well in our data-driven causal measure.

\begin{table}[h]
\begin{tabular}{|l|l|l|}
\hline
 & $ \Omega_2 \ncauses \Omega_3$ & $ \Omega_3 \ncauses \Omega_2$ \\ 
 \hline
\textit{Symmetric} & 0.90  & 0.97  \\ 
\hline
\textit{Asymmetric} & 0.00 & 0.05  \\ 
\hline
\textit{Asymmetric w/ Correction} & 1.00 & 0.89 \\ 
\hline
\end{tabular}
\caption{Evaluation of our Koopman causality measure on the \rossler confounder model. Koopman causality is evaluated directly for the \textit{Symmetric} and \textit{Asymmetric} cases, while in the \textit{Asymmetric w/ Correction} case we adjust for the $\M_1$ confounder.\label{tab:confouder}}
\end{table}

\subsection{Asymmetric coupling}
Success in the symmetric case shows the promise of our method. However, our current approximation algorithm is relatively simple, and as we now show it is not generically immune to spurious association from confounding. Switching the \rossler confounder to asymmetric coupling, with $c_2 = 0.5$ and $c_3 = 0.3$, the middle row of Table~\ref{tab:confouder} shows that our Koopman causality measure erroneously identifies a causal relation between $\M_2$ and $\M_3$ in most cases. 

\subsection{Correcting for Confounding}
Given that our theoretical framework is indifferent to confounding, we are optimistic that improved approximation methods will be able to more reliably distinguish spurious non-causal associations. For now though, our current method can correct for spurious associations by adjusting for confounding, given that the confounders are observed. 

To adjust for confounding, we utilize the multi-variate ability of our Koopman causality method. Consider a causal query of interest, $\M_C \causes \M_E$ in the presence of a confounding component $\M_\xi$, for example the causal relation between $\M_2$ and $\M_3$ (either direction) in the \rossler confounder with confounder $\M_\xi = \M_1$. If the initial Koopman causality measure identifies a causal relation, $\M_C \ddcauses \M_E$, but we suspect there is no actual causal relation, we can check whether the detected value is spurious by comparing it to other causal measures involving the confounder. 

\begin{figure}
    \begin{center}
        \includegraphics[width=0.48 \textwidth]{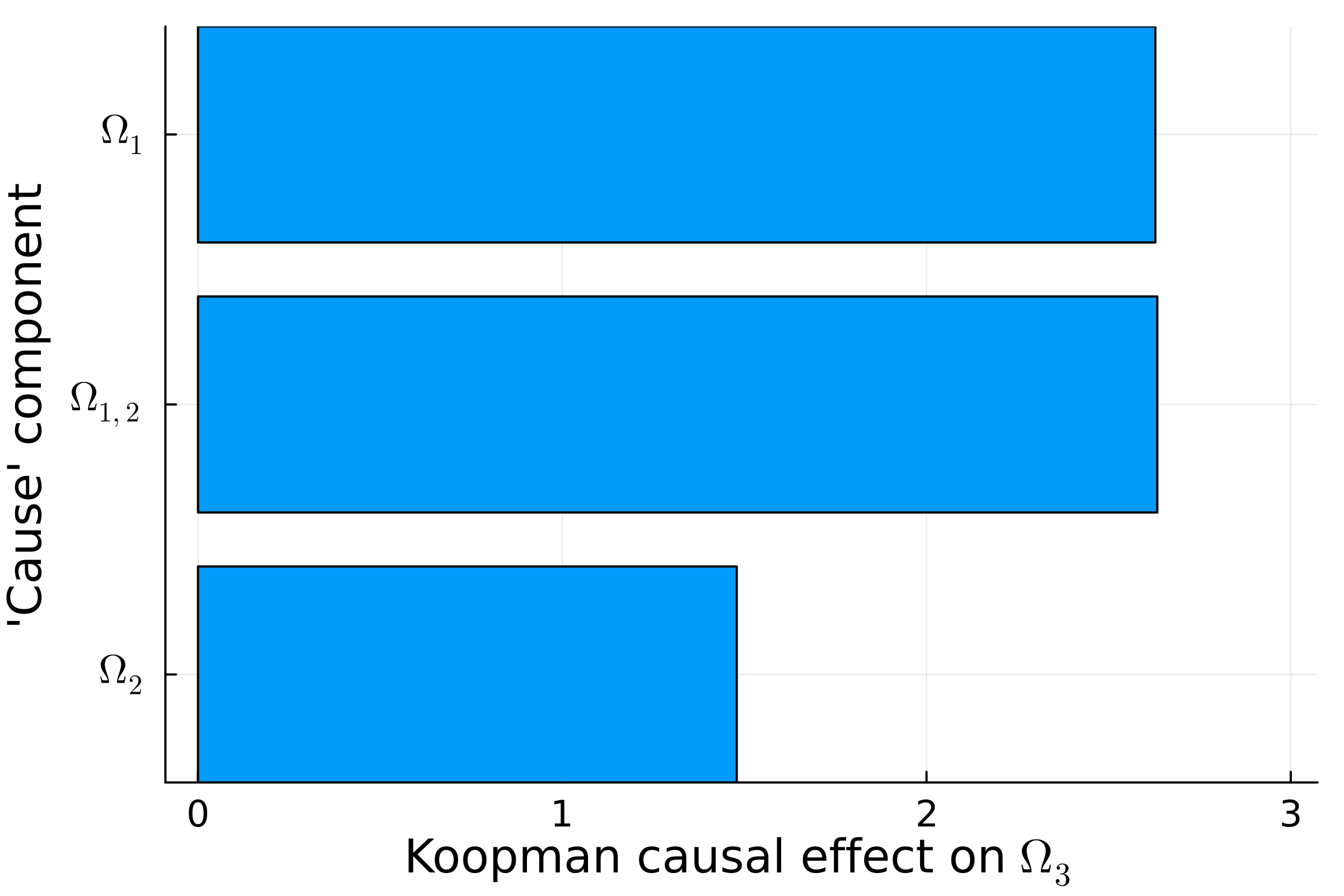}
    \end{center}
    \caption{Correcting for confounding in the asymmetric \rossler model.}
    \label{fig:confounding_correction}
\end{figure}

In particular, we compute three causal measures: $\M_C \ddcauses \M_E$, $\M_\xi \ddcauses \M_E$, and $\M_{C, \xi} \ddcauses \M_E$. 
An example from the \rossler confounder model is shown in Figure~\ref{fig:confounding_correction} for $\M_2 \causes \M_3$ with confounding from $\M_\xi = \M_1$. We see that, while $\M_2 \ddcauses \M_3$ is significantly non-zero, $\M_1 \ddcauses \M_3 = \M_{1,2} \ddcauses \M_3$. This tells us that $\M_1 \causes \M_3$ is the true causal relation, and the signal we see in $\M_2 \ddcauses \M_3$ is a residual effect from the spurious association due to the confounding from $\M_1$. This analysis is reminiscent of the \emph{backdoor adjustment} from the statistical approach to causal inference~\cite{pearl09a}, in which confounders are accounted for by conditioning on them.

We can understand how our correction works by using the lens of the \emph{partial information decomposition}~\cite{will10a}, which formalizes the intuitive notions of \emph{unique}, \emph{redundant}, and \emph{synergistic} effects of two components on a third. First, consider the measured effect of component $\M_2$ in Figure~\ref{fig:confounding_correction}. When we then also add $\M_1$ we see a significantly larger effect for the combined $\M_{1,2}$ component. Now, consider starting with the measured effect of $\M_1$. When we then add $\M_2$, we see the combined effect of $\M_{1,2}$ is the same as that of $\M_1$ alone. Therefore, we see that $\M_1$ contributes a \emph{unique} causal effect on $\M_3$ and any measured effect from $\M_2$ is \emph{redundant} with the effect from $\M_1$. And, adding $\M_2$ contributes nothing beyond what is measured from the effect of $\M_1$, therefore there is no synergistic effect from $\M_1$ and $\M_2$ combined. We thus conclude that $\M_2$ carries no unique causal effect on $\M_3$, and the measured value $\M_2 \ddcauses \M_3$ is spurious. 

To turn this observation into a concrete correction to determine whether a measured causal effect is spurious, we do the following. First, we compute the three causal Koopman measures $\M_C \ddcauses \M_E$, $\M_\xi \ddcauses \M_E$, and $\M_{C,\xi} \ddcauses \M_E$. Since we expect that $\M_\xi$ is a confounder, and thus has a positive causal effect, we use $\M_\xi \ddcauses \M_E$ as a baseline to compare against. Define,
\begin{align*}
    \Delta_{C} &:= |(\M_C \ddcauses \M_E) - (\M_{\xi} \ddcauses \M_E)|\\
    \Delta_{C,\xi} &:= | (\M_{C,\xi} \ddcauses \M_E) - (\M_\xi \ddcauses \M_E)|
    ~.
\end{align*}
Given the above discussion, if $\M_C \ddcauses \M_E$ is spurious due to confounding, we expect $\Delta_C >> \Delta_{C,\xi}$ since $\M_{C,\xi} \ddcauses \M_E$ should be very similar to $\M_\xi \ddcauses \M_E$ but $\M_C \ddcauses \M_E$ should be significantly smaller than $\M_\xi \ddcauses \M_E$. Therefore, we set an order-of-magnitude threshold to determine whether the measured value $\M_C \ddcauses \M_E$ is spurious due to confounder $\M_\xi$:
\begin{align*}
    \M_C \ncauses \M_E \text{ if } \Delta_C \geq 10 \; \Delta_{C,\xi}
    ~.
\end{align*}
Using this correction criterion on the asymmetric \rossler confounder tests, the bottom row of Table \ref{tab:confouder} shows that we can now correctly identify that there is no actual causal relationship in most cases.

\section{Application to Spatiotemporal Systems}
\label{sec:examples}

Up to this point we have used the coupled \rossler system as the guiding example to demonstrate our theory and numerical causal measure. We now explore more complex systems to highlight the utility of our Koopman-based causal measure. In particular, we use the Lorenz 96 system as an idealized spatiotemporal system to show how our measure can discover structural properties of a system directly from data. In particular, asymmetrical coupling between oscillators in the Lorenz 96 model leads to an asymmetrical causal flow over space through time that our causal measure can identify. 
We will see that our Koopman causality measure identifies a causal wave that propagates in the same direction and at roughly the same speed that information propagates in the system, as identified through perturbation propagation. 

\subsection{The Lorenz 96 System}

\begin{figure}
    \begin{center}
        \includegraphics[width=0.5 \textwidth, trim={0cm 4.0cm 9.5cm 0.0cm}, clip]{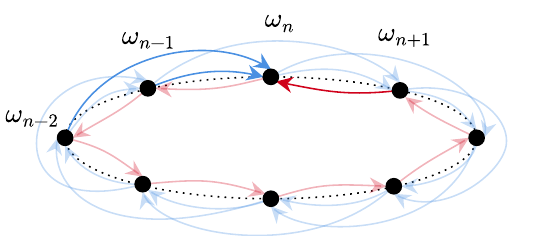}
    \end{center}
    \caption{Asymmetry in the spatial coupling between oscillators of the Lorenz 96 model.}
    \label{fig:L96}
\end{figure}

The Lorenz 96 system is an idealized model of the dynamics of a scalar quantity in the atmosphere~\cite{Keri22a}. It consists of $N$ sites on a ring, which are viewed as $N$ equally spaced locations on a single line of latitude on the Earth. Each degree of freedom $\omega_n$ on an individual site $n$ evolve according to 
\begin{align}
    \dot{\omega}_n = (\omega_{n+1} - \omega_{n-2})\omega_{n-1} - \omega_n + F
    ~.
    \label{eq:L96}
\end{align}
We will refer to the $\omega_n$ simply as `oscillators'. The last term $F \geq 0$ is the \emph{forcing}, as it is the rate at which the scalar quantity is added uniformly into the system. It is the only free parameter of this formulation of the model and controls how chaotic the system orbits are. The middle term, $-\omega_n$, is the \emph{dissipation} term that dictates the rate of loss of the scalar quantity over time. 


\begin{figure}
    \begin{center}
        \includegraphics[width=0.5 \textwidth]{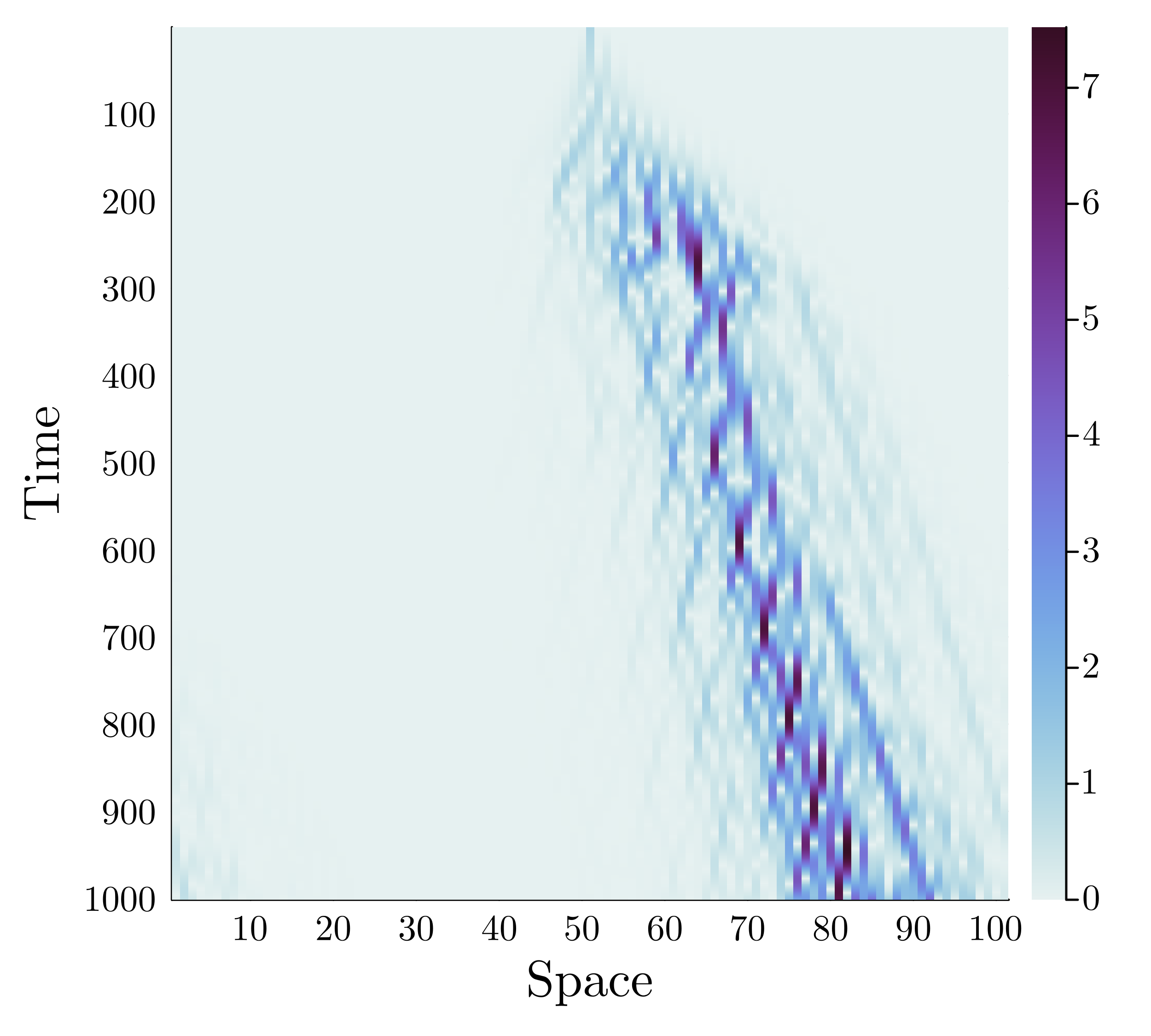}
    \end{center}
    \caption{Perturbation propagation in the Lorenz 96 model. The spacetime diagram depicts the (absolute value) difference between two orbits of the Lorenz 96 model that have identical initial conditions except for a perturbation in the central oscillator. This experiment indicates that information propagates in the system as a traveling wave in the clockwise direction.}
    \label{fig:L96_pert}
\end{figure}

Finally, the most important term for our purposes here is the first, $(\omega_{n+1} - \omega_{n-2})\omega_{n-1}$. Because of the asymmetry in spatial coupling, this is the \emph{advection} term. The dynamics of the oscillator $\omega_n$ at site $n$ depends only on the nearest neighbor in the clockwise direction $\omega_{n+1}$, but it depends on both the nearest $\omega_{n-1}$ and next-nearest $\omega_{n-2}$ neighbors in the counterclockwise direction; see Figure~\ref{fig:L96}. The imbalance in spatial coupling leads to a clockwise transport of the quantity in the system over time. 

For our experiments we use $N = 101$ oscillators and a moderate forcing value of $F = 4.0$. 

Due to the simplicity of the model, intuition about material transport translates closely to a `causal flow' in the Lorenz 96 system. For now, consider each individual oscillator as its own component, i.e. $\M_n = \{\omega_n\}$. As time increases, the component flow maps $P_n \Phi^t$ will depend on more components in the counterclockwise direction than the clockwise due to the asymmetric coupling discussed above. Thus, there is a dominant clockwise causal flow in the system similar to the clockwise material transport. 

\subsubsection{Perturbation Propagation and Information Flow}

Before applying our Koopman causality measure, we first identify a baseline for comparison. We perform a perturbation propagation experiment in which we evolve a Lorenz 96 model from random initial conditions and then evolve a second model with the same initial conditions except for a perturbation applied to the oscillator at the middle site $n = 51$. As discussed by, e.g., Ref.~\cite{pack85a}, perturbation spreading indicates how \emph{information} propagates through space over time in the system. A more recent example is the popular ``slinky drop'' experiment~\cite{cross12a,slinky24a}. A long slinky is held vertical, from the top, and then dropped. Slow motion video reveals that the bottom of the slinky does not move until the rest of the slinky reaches it. A common explanation is that the bottom of the slinky ``does not know'' that the top has been released until that information reaches the bottom, and information propagates at the (relatively slow) longitudinal wave speed of the slinky. 


A spacetime diagram of perturbation spreading in the Lorenz 96 model is shown in Figure~\ref{fig:L96_pert}, with spatial sites on the horizontal axis and time on the vertical increasing from top to bottom. From this, we see that information travels almost exclusively in the clockwise direction (increasing $n$), at a finite speed. Similar to the slinky drop, the other oscillators do not instantly ``know'' that the central oscillator has been perturbed. Information about the central oscillator perturbation travels to distant oscillators at a finite rate, due to the locality of interactions between the oscillators. It is notable that the local coupling occurs in \emph{both} directions, as depicted in Figure~\ref{fig:L96}, but the spatial asymmetry and cyclic feedback of the local coupling conspire to give a net flow of information in the clockwise direction, as seen in Figure~\ref{fig:L96_pert}.




\begin{figure*}
    \begin{center}
        \includegraphics[width=0.8 \textwidth]{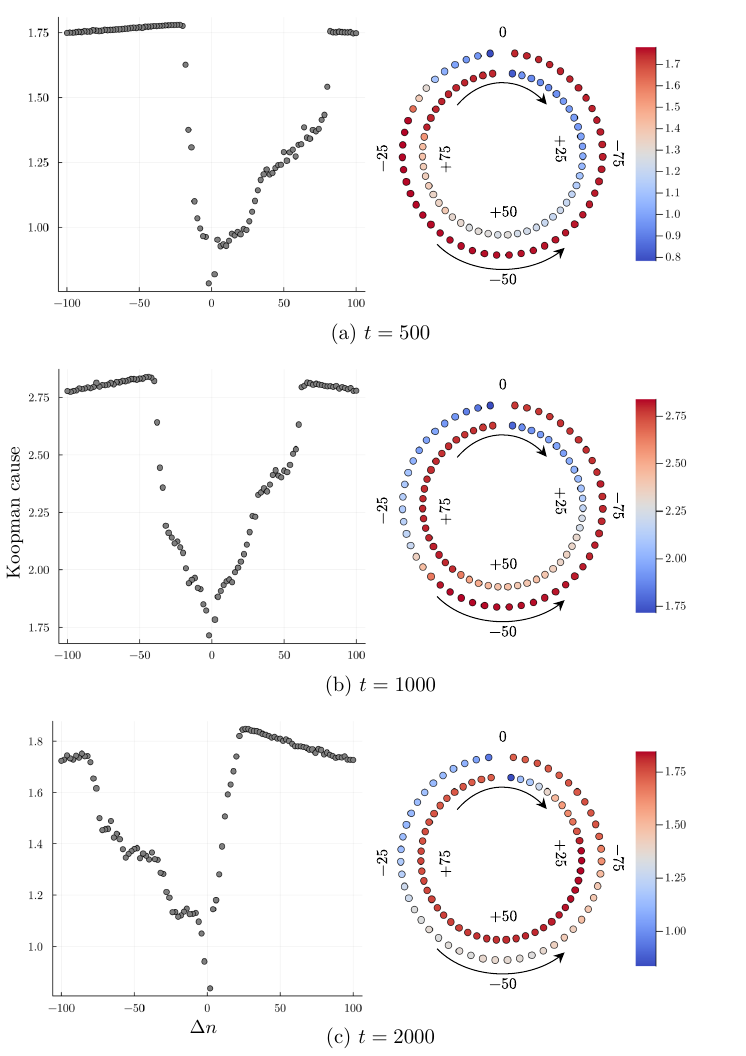}
    \end{center}
    \caption{Cumulative causal effects in the Lorenz 96 model.}
    \label{fig:L96_cause}
\end{figure*}

\subsubsection{Koopman Analysis and Causal Flow}

The information flow analysis just described is achieved through perturbing the system. Similar to counterfactual experiments, this requires exact control over the system to very high precision, which is not realistic in practice. We now apply our Koopman causal measure to the Lorenz 96 system, and again emphasize it is purely observational and does not require interventions (e.g., perturbations). Information flow is often assumed to be synonymous with causal flow such that causality is identified through measures of information flow. For example, transfer entropy, a measure of information flow, is a nonlinear generalization of Wiener-Granger causality. Through comparison with the perturbation propagation in Figure~\ref{fig:L96_pert}, our Koopman causality measure reveals that there is indeed a close correspondence between causal flow and information flow in the Lorenz 96 system. 


Consider a `target site' $n^*$---it does not matter which specific site we choose since the Lorenz 96 system, as a whole, is rotationally symmetric. We want to identify the causal effects of other oscillators on the target oscillator $\omega_{n^*}$, so we designate the `effect' component $\M_E$ as the singelton $\{\omega_{n^*}\}$. 
Because causal flows must follow the spatial ring structure of the Lorenz 96 model, oscillators on the opposite side of the ring do not ``directly'' influence the target oscillator through their local interactions. Their causal influence must propagate over space to reach the target oscillator indirectly through the oscillators between them. This means there is \textbf{significant confounding} between oscillators. However, given that all oscillators are observed and that we know the connection topology if Figure~\ref{fig:L96}, we can account for the confounding and unravel how the intertwined nonlinear relationships unfold over time. 

We do so by examining the \emph{cumulative} causal effects of increasingly-large sets of neighbors on the target oscillator. For the `cause' component, we take all neighbors to and including the \emph{neighbor number} $\Delta n$, with the direction of neighbors taken given by the sign of the neighbor number. For example, with $\Delta n = 2$ we take the clockwise nearest and next-nearest neighbors for the `cause' component $\M_C$. For $\Delta n = -4$ we have $\M_C = \{\omega_{n^* -1}, \omega_{n^* -2}, \omega_{n^* - 3}, \omega_{n^* - 4}\}$, etc. 

Figure~\ref{fig:L96_cause} shows results of $\M_C \ddcauses \M_E$ for the cumulative $\M_C$ components at three different time shifts, $t = \{500, 1000, 2000\}$. For each time shift (each row in the Figure), there are two different plots that display the same results. The scatter plots in the left column show the causal effect values $\M_C \ddcauses \M_E$ on the vertical axis, and the horizontal axis is the neighbor number $\Delta n$ that gives cumulative $\M_C$ component just defined. Thus, moving left from the center of each plot shows the cumulative causal effect of increasingly many counterclockwise neighbors out to $\Delta n = -100$, for which $\M_C$ contains all other oscillators aside from the target. Similarly, moving right from center shows the cumulative causal effect of increasingly many clockwise neighbors out to $\Delta n = 100$, for which $\M_C$ again contains all other oscillators aside from the target. 

The circle plots in the right column also show $\M_C \ddcauses \M_E$ for the same values of $\Delta n$. These are spatial embeddings of the scatter plots, with the value of the causal effect given by the color on each plot marker. The outer ring is equivalent to the left side of the scatter plots with negative $\Delta n$---the counterclockwise neighbors. The inside ring is the clockwise neighbors with positive $\Delta n$, and is equivalent to the right half of the scatter plots. 

Collectively in Figure~\ref{fig:L96_cause} we see a plateau in the value of the cumulative causal effect $\M_C \ddcauses \M_E$ that travels counterclockwise over time, which corresponds to a propagating front of \emph{causal flow} that travels clockwise over time. This is easiest to see from inspecting the cumulative causal effect of the counterclockwise neighbors---those with negative neighbor numbers. Again, we are measuring the cumulative effect of neighbors on the target oscillator, and counterclockwise neighbors to the target is in the clockwise direction. So the cumulative effect of counterclockwise neighbors is measuring in the same direction that the front is traveling.  

Starting at the shortest time step $t = 500$, the causal flow front has the least amount of time to travel and so travels the shortest distance. Therefore $t = 500$ has the fewest counterclockwise neighbors that can effect the target through the propagating front of causal flow. Any additional counterclockwise neighbors added cannot contribute via the front of causal flow, hence the plateau. As we increase to $t = 1000$ and $t = 2000$ we see a smaller and smaller plateau because the front has more time to travel further, and hence more counterclockwise neighbors can effect the target oscillator via the front. This is exactly what we see in Figure~\ref{fig:L96_cause}. Significantly, the spatial location of the front at $t = 500$ (around $\Delta n = -20$) and $t = 1000$ (around $\Delta n = -30$) agree with perturbation front locations in Figure~\ref{fig:L96_pert}. 

The cumulative causal effect of clockwise neighbors measures causation that does not follow the clockwise causal flow front, at least for small enough neighbor numbers. Because the Lorenz 96 model is on a ring, if you add enough clockwise neighbors you eventually come to the causal flow front from behind, where we again see a plateau in the values of our causal measure $\M_C \ddcauses \M_E$. The plateau in the clockwise direction occurs at roughly the same location on the ring as in the counterclockwise direction, as best seen from the circle plots in the right column of Figure~\ref{fig:L96_cause}. 

\begin{figure}
    \begin{center}
        \includegraphics[width=0.5 \textwidth]{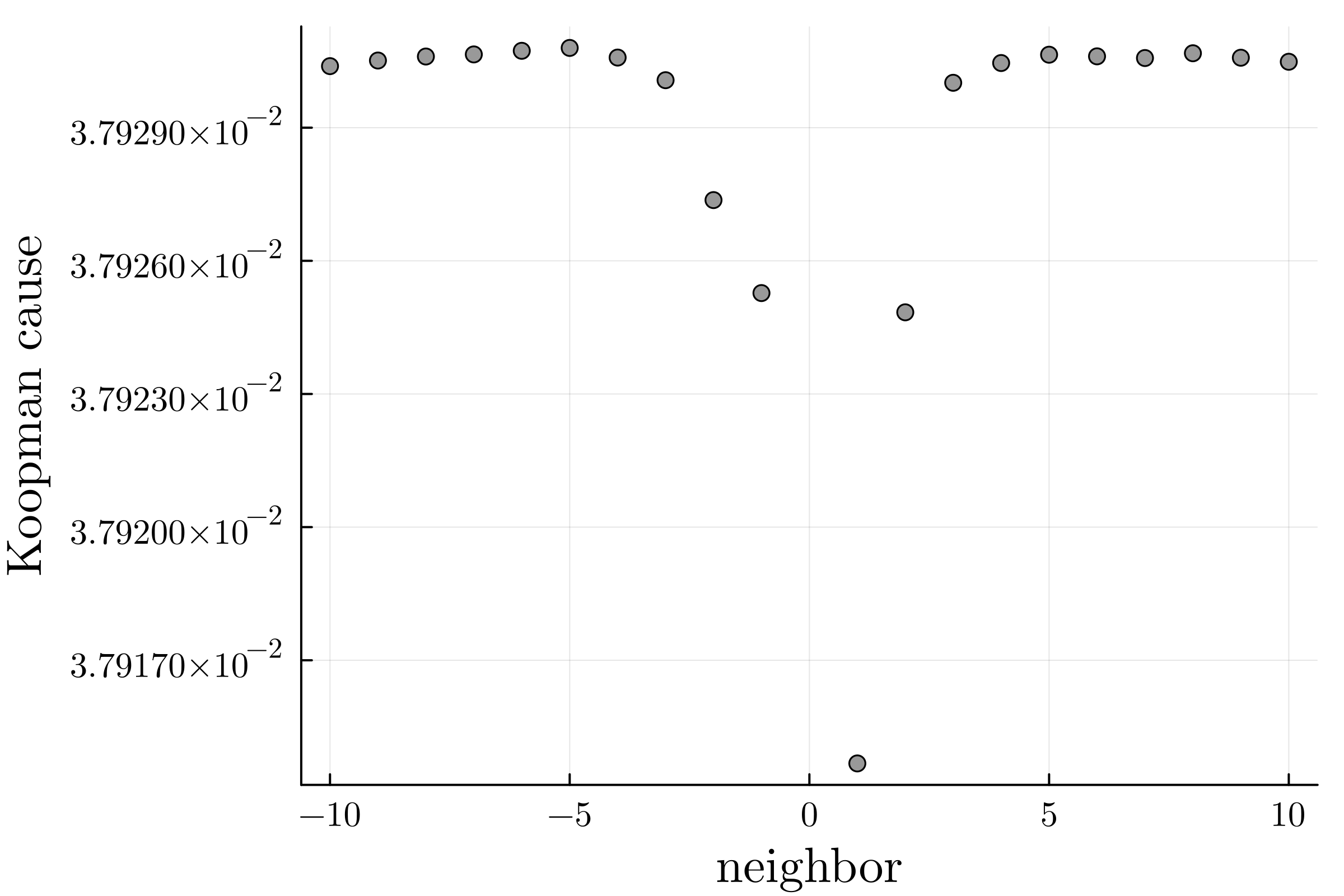}
    \end{center}
    \caption{Cumulative causal effects approximating the differential time scale.}
    \label{fig:L96_instant}
\end{figure}

To summarize, Figure~\ref{fig:L96_cause} shows that our Koopman causality measure identifies the dominant clockwise causal flow in the Lorenz 96 system. And, the causal flow we identify travels in the same direction and at roughly the same speed as the flow of information identified by perturbation propagation in Figure~\ref{fig:L96_pert}.


\subsubsection{Instantaneous Causal Effects and Direct Influence}

The above analysis implicitly assumed knowledge of the spatial ring structure of the Lorenz 96 model. This structure arises from the coupling between oscillators in the advection term of the equations of motion (\ref{eq:L96}). Recall that for continuous-time dynamical systems, the differential equations of motion are the time derivatives of the flow maps. And, since causality in dynamical systems arises from the flow maps, the differential equations of motion represent \emph{instantaneous causal effects}. There is often an emphasis in causality literature on `direct' vs `indirect' causal influence \cite{sun15a,rung19a}. From our rigorous formulation of dynamical system causality, it is clear that the use of `direct' causation in dynamical systems really means `instantaneous' causation. `Indirect' causation occurs at later times when the flow maps pickup additional dependencies not present in the differential equations of motion. 

To investigate instantaneous causation using our Koopman causality measure, we perform the same cumulative causal analysis as above at a very small time scale---we use an equivalent of $t=0.1$ (i.e., one time step using an integration time of $\delta = 0.001$ compared to $\delta = 0.01$ used in Figure~\ref{fig:L96_cause}). The results are shown in Figure.~\ref{fig:L96_instant}. We can see that there is an asymmetry between clockwise and counterclockwise and the causal effect plateaus quickly in both directions. However, it is not clear from this analysis that the instantaneous causation includes one clockwise neighbor and two counterclockwise neighbors. 

Because instantaneous causal effects arise from the differential equations of motion, they can be identified using equation discovery methods \cite{Crut87a,brun16b}. Moreover, the time derivative of the Koopman semigroup is the Koopman generator---the Koopman equivalent of the equations of motion. It has been shown that equation discovery algorithms are special cases of approximation methods for the Koopman generator \cite{klus20a}. Thus, in principle, our Koopman method is capable of identifying direct causal effects in some cases. The results in Figure~\ref{fig:L96_instant} are essentially equivalent to approximating the Koopman generator from finite differences, using marginal and joint dictionaries. We thus find that this straightforward method is not sufficient for approximating the Koopman generator of this Lorenz 96 model, even using very small differences. Other methods, such as automatic differentiation, could potentially find the correct instantaneous causal relations. 

To close, remember that in many cases of interest we \emph{do} know the governing equations and thus the topology of direct causal relations (though not necessarily their magnitudes). What is not clear, for many complex systems, is how these causal relations unfold as time evolves. We again emphasize that the flow maps for these systems are complicated and do not have closed-form analytic solutions. The strength of our data-driven Koopman causality method is its ability to discover causal relations at non-instantaneous time scales.

\section{Discussion}
\label{sec:discussion} 

The existing tools for causality in dynamical systems largely fall into two categories based on either i) information flows, or ii) delay coordinate embeddings. Having introduced our Koopman theory of causality, we now overview how it relates to these other methods. Conceptually, it is most similar to information flows. In fact, information flows in deterministic dynamical systems can be formulated in terms of Perron-Frobenius operators \cite{liang05a}, which are dual to Koopman operators \cite{laso94a}---they evolve probability distributions over systems states, rather than observables of the system states. 

Note that information flow is a more general concept that extends beyond causal discovery, and so measures to quantify information flow, such as transfer entropy~\cite{schrei00a}, are not necessarily designed to be measures of causality. Additional axioms are typically posited to identify causal flows through information flows, notably \emph{asymmetry}---the causal flow from $\M_2$ to $\M_1$ is not necessarily the same as from $\M_1$ to $\M_2$---and \emph{nil causality}---if $\M_2 \ncauses \M_1$ then the causal flow from $\M_2$ to $\M_1$ should be zero. From Figure~\ref{fig:causalDMD} we can see that our Koopman causality measure satisfies these axioms. Additional work is necessary to more rigorously connect our Koopman causality framework to information flows through the duality of Koopman and Perron-Frobenius operators. 

The other category of causal methods for dynamical systems are those based on delay coordinate embeddings, most notably convergent cross mapping~\cite{sugi12a}. Whereas information flows are more of a `top down' approach to causality, delay embedding methods are `bottom up', with its definition of causality rigorously based on the ability of past histories of some degrees of freedom to encode the effects of other degrees of freedom in the system \cite{Take81}. Interestingly, while our Koopman theory is similarly `bottom up', starting from the rigorous Definition~\ref{def:dyn_cause} of dynamical system causality, it does \emph{not} utilize delay embeddings. 
Because delay embeddings of a given observable can encode the \emph{entirety} of the remaining degrees of freedom \cite{rupe22b}, they cannot identify the unique causal influence of a single component in the presence of additional unobserved causal components (a non-empty `remainder'). Our approach utilizing Koopman operators is advantageous in this respect for high-dimensional systems with many interacting components. 

Finally, we note that the practical application of our data-driven causality measure is limited by the quality of the DMD approximation to the Koopman operators used in the marginal and joint models. If the DMD models give very poor approximations to the evolution of the identity observable (or other test functions in $\F_E$, more generally), then the comparison of marginal and joint model errors is not necessarily reflective of causal relations in the system. In particular, the class of DMD models used in this work do not reliably give good approximations for chaotic systems. 

That said, our framework allows for the use of more sophisticated Koopman approximations to give more robust data-driven measures of causality. Alternatively, it may be possible to use our Koopman theory to define a measure of causality based on defects between component subspaces. If $\M_C \ncauses \M_E$ then $\Koop^t \F_E \subseteq \F_{E,R}$, and so a possible route to measure causal effects would be to identify the defect between $\Koop^t \F_E$ and $\F_{E,R}$. The ``principal angle'' between these space~\cite{deutsch1995angle} is a promising approach to quantify this defect from data. 

\section{Conclusion}
\label{sec:conclusion}

We have introduced a new theory of causality in nonlinear dynamical systems using Koopman operators. Our theory is based on a rigorous definition of causal mechanism in dynamical systems using flow maps, analogous to the standard definition of causal mechanism in structural causal models. We defined causality in the Koopman framework as flows between function subspaces, and proved that this is equivalent to the definition of causal mechanism using flow maps. Our Koopman theory of causality lends itself to a straightforward data-driven measure of causality based on the dynamic mode decomposition algorithm. After giving some basic demonstrations on the coupled \rossler system, we showed that our data-driven causality measure can identify causal flow in the Lorenz 96 model associated with advective transport, in agreement with information flow identified in a perturbation propagation experiment.

In addition to providing a rigorous foundation of causal mechanisms for dynamical systems, our data-driven Koopman causality measure provides new capabilities not found in most causality algorithms for dynamical systems. The definition of causal mechanism in terms of flow maps clarifies the importance of time scale, and our Koopman method identifies relations at a chosen time scale. It also naturally incorporates multivariate (or polyadic) causal relationships \cite{james16a} through the flexible definition of components. Because DMD can scale to high-dimensional systems, our method can evaluate causal relations between high-dimensional components. Time series causality methods, by contrast, first require a dimensionality reduction step, and so can only evaluate causal relationships between the reduced components. 

\section*{Acknowledgments}

This work was supported by the U.S. Department of Energy (DOE), Office of Science, Office of Biological and Environmental Research, Regional and Global Model Analysis program area as part of the HiLAT-RASM project.
The Pacific Northwest National Laboratory (PNNL) is operated for DOE by Battelle Memorial Institute under contract DE-AC05-76RLO1830.

\bibliography{ergodic}

@inproceedings{seabold2010statsmodels,
  title={statsmodels: Econometric and statistical modeling with python},
  author={Seabold, Skipper and Perktold, Josef},
  booktitle={9th Python in Science Conference},
  year={2010},
}

@article{sugihara2012detecting,
  title={Detecting causality in complex ecosystems},
  author={Sugihara, George and May, Robert and Ye, Hao and Hsieh, Chih-hao and Deyle, Ethan and Fogarty, Michael and Munch, Stephan},
  journal={science},
  volume={338},
  number={6106},
  pages={496--500},
  year={2012},
  publisher={American Association for the Advancement of Science}
}

@article{will10a,
  title={Nonnegative decomposition of multivariate information},
  author={Williams, P. L. and Beer, R. D.},
  journal={arXiv preprint arXiv:1004.2515},
  year={2010}
}

@article{Julia17,
    title={Julia: A fresh approach to numerical computing},
    author={Bezanson, J. and Edelman, A. and Karpinski, S. and Shah, V.B.},
    journal={SIAM {R}eview},
    volume={59},
    number={1},
    pages={65--98},
    year={2017},
    publisher={SIAM},
    doi={10.1137/141000671},
    url={https://epubs.siam.org/doi/10.1137/141000671}
}

@article{cross12a,
  title={Modeling a falling slinky},
  author={Cross, R.C. and Wheatland, M.S.},
  journal={American journal of physics},
  volume={80},
  number={12},
  pages={1051--1060},
  year={2012},
  publisher={AIP Publishing}
}

@inbook{bare22a,
author = {Bareinboim, E. and Correa, J. D. and Ibeling, D. and Icard, T.},
title = {On {P}earl’s Hierarchy and the Foundations of Causal Inference},
year = {2022},
isbn = {9781450395861},
publisher = {Association for Computing Machinery},
address = {New York, NY, USA},
edition = {1},
booktitle = {Probabilistic and Causal Inference: The Works of Judea Pearl},
pages = {507–556},
numpages = {50}
}

@misc{slinky24a,
  title = {Slinky Drop Answer, {Veritasium}},
  author = {D. Muller},
  howpublished = {\url{https://www.youtube.com/watch?v=eCMmmEEyOO0}},
  year = {Accessed: 2024-10-02}
}

@article{pack85a,
  title={Two-dimensional cellular automata},
  author={Packard, N. H. and Wolfram, S.},
  journal={Journal of Statistical physics},
  volume={38},
  number={5},
  pages={901--946},
  year={1985},
  publisher={Springer}
}

@incollection{deutsch1995angle,
  title={The angle between subspaces of a Hilbert space},
  author={Deutsch, F.},
  booktitle={Approximation theory, wavelets and applications},
  pages={107--130},
  year={1995},
  publisher={Springer}
}

@article{sun15a,
  title={Causal network inference by optimal causation entropy},
  author={Sun, J. and Taylor, D. and Bollt, E. M.},
  journal={SIAM Journal on Applied Dynamical Systems},
  volume={14},
  number={1},
  pages={73--106},
  year={2015},
  publisher={SIAM}
}

@article{hernan08a,
  title={Does obesity shorten life? The importance of well-defined interventions to answer causal questions},
  author={Hern{\'a}n, M. A. and Taubman, S. L.},
  journal={International journal of obesity},
  volume={32},
  number={3},
  pages={S8--S14},
  year={2008},
  publisher={Nature Publishing Group}
}

@article{schrei00a,
  title={Measuring information transfer},
  author={Schreiber, T.},
  journal={Physical review letters},
  volume={85},
  number={2},
  pages={461},
  year={2000},
  publisher={APS}
}

@article{tsoni18a,
  title={Convergent cross mapping: theory and an example},
  author={Tsonis, A. A. and Deyle, E. R. and Ye, H. and Sugihara, G.},
  journal={Advances in nonlinear geosciences},
  pages={587--600},
  year={2018},
  publisher={Springer}
}

@article{liang05a,
  title={Information transfer between dynamical system components},
  author={Liang, X. S. and Kleeman, R.},
  journal={Physical review letters},
  volume={95},
  number={24},
  pages={244101},
  year={2005},
  publisher={APS}
}

@article{james16a,
  title={Information flows? A critique of transfer entropies},
  author={James, R. G. and Barnett, N. and Crutchfield, J. P.},
  journal={Physical review letters},
  volume={116},
  number={23},
  pages={238701},
  year={2016},
  publisher={APS}
}

@article{ay08,
  title={Information flows in causal networks},
  author={Ay, N. and Polani, D.},
  journal={Advances in complex systems},
  volume={11},
  number={01},
  pages={17--41},
  year={2008},
  publisher={World Scientific}
}

@article{roe09a,
  title={Feedbacks, timescales, and seeing red},
  author={Roe, G.},
  journal={Annual Review of Earth and Planetary Sciences},
  volume={37},
  number={1},
  pages={93--115},
  year={2009},
  publisher={Annual Reviews}
}

@book{pearl09a,
  title={Causality},
  author={Pearl, J.},
  year={2009},
  publisher={Cambridge university press}
}

@article{camp23a,
  title={Discovering causal relations and equations from data},
  author={Camps-Valls, G. and Gerhardus, A. and Ninad, U. and Varando, G. and Martius, G. and Balaguer-Ballester, E. and Vinuesa, R. and Diaz, E. and Zanna, L. and Runge, J.},
  journal={Physics Reports},
  volume={1044},
  pages={1--68},
  year={2023},
  publisher={Elsevier}
}

@article{rupe24a,
  title={On principles of emergent organization},
  author={Rupe, A. and Crutchfield, J. P.},
  journal={Physics Reports},
  volume={1071},
  pages={1--47},
  year={2024},
  publisher={Elsevier}
}

@article{gonz21a,
  title={The kernel perspective on dynamic mode decomposition},
  author={Gonzalez, E. and Abudia, M. and Jury, M. and Kamalapurkar, R. and Rosenfeld, J. A.},
  journal={arXiv preprint arXiv:2106.00106},
  year={2021}
}

@article{wien56a,
  title={The theory of prediction},
  author={N., Wiener},
  journal={Modern Mathematics for Engineers; Beckenbach, E.F., Ed.; McGraw-Hill: New York, NY, USA},
  year={1956}
}

@article{gran69,
  title={Investigating causal relations by econometric models and cross-spectral methods},
  author={Granger, C. W.J.},
  journal={Econometrica: journal of the Econometric Society},
  pages={424--438},
  year={1969},
  publisher={JSTOR}
}

@article{Harn17a,
  title = {Topological Causality in Dynamical Systems},
  author = {Harnack, D. and Laminski, E. and Sch\"unemann, M. and Pawelzik, K. R.},
  journal = {Phys. Rev. Lett.},
  volume = {119},
  issue = {9},
  pages = {098301},
  numpages = {5},
  year = {2017},
  month = {Sep},
  publisher = {American Physical Society},
  doi = {10.1103/PhysRevLett.119.098301}
}

@article{merl14a,
  title={Interacting components of the top-of-atmosphere energy balance affect changes in regional surface temperature},
  author={Merlis, T. M.},
  journal={Geophysical Research Letters},
  volume={41},
  number={20},
  pages={7291--7297},
  year={2014},
  publisher={Wiley}
}

@misc{Gund19a,
  title = {Random {Fourier} Features},
  howpublished = {\url{https://gregorygundersen.com/blog/2019/12/23/random-fourier-features/}},
  author = {G. Gundersen}, 
  year = {2019},
  note = {Accessed: 2010-09-30}
}

@book{pauls16a,
  title={An introduction to the theory of reproducing kernel {H}ilbert spaces},
  author={Paulsen, V. I. and Raghupathi, M.},
  volume={152},
  year={2016},
  publisher={Cambridge university press}
}

@article{Keri22a,
title = {On the {L}orenz '96 model and some generalizations},
journal = {Discrete and Continuous Dynamical Systems - B},
volume = {27},
number = {2},
pages = {769-797},
year = {2022},
issn = {1531-3492},
doi = {10.3934/dcdsb.2021064},
author = {J. Kerin and H. Engler}
}

@book{rudi94,
  title={Fourier analysis on groups},
  author={Rudin, W.},
  year={1994},
  publisher={Wiley-Interscience},
}

@article{rahi07a,
  title={Random features for large-scale kernel machines},
  author={Rahimi, A. and Recht, B.},
  journal={Advances in neural information processing systems},
  volume={20},
  year={2007}
}

@article{Koop31a,
  title={Hamiltonian systems and transformation in {H}ilbert space},
  author={Koopman, B. O.},
  journal={Proc. Natl. Acad. Sci. USA},
  volume={17},
  number={5},
  pages={315},
  year={1931},
  publisher={National Academy of Sciences}
}

@article{brun22a,
  title={Modern {K}oopman Theory for Dynamical Systems},
  author={Brunton, S. L. and Budi\v{s}i\'{c}, M. and Kaiser, E. and Kutz, J. N.},
  journal={SIAM Review},
  volume={64},
  number={2},
  pages={229--340},
  year={2022},
  publisher={Society for Industrial and Applied Mathematics}
}

@article{rupe22b,
doi = {10.1088/1367-2630/ac95b7},
year = {2022},
publisher = {IOP Publishing},
volume = {24},
number = {10},
pages = {103033},
author = {A. Rupe and V. V. Vesselinov and J. P. Crutchfield},
title = {Nonequilibrium statistical mechanics and optimal prediction of partially-observed complex systems},
journal = {New Journal of Physics},
}

@book{laso94a,
  title={Chaos, fractals, and noise: stochastic aspects of dynamics},
  author={Lasota, A. and Mackey, M. C.},
  volume={97},
  year={1994},
  publisher={Springer Science}
}

@article{will15a,
  title={A data--driven approximation of the {K}oopman operator: Extending dynamic mode decomposition},
  author={Williams, M. O. and Kevrekidis, I. G. and Rowley, C. W.},
  journal={Journal of Nonlinear Science},
  volume={25},
  number={6},
  pages={1307--1346},
  year={2015},
  publisher={Springer}
}

@article{klus16a,
  title={On the numerical approximation of the Perron-Frobenius and Koopman operator},
  author={Klus, S. and Koltai, P. and Sch{\"u}tte, C.},
  journal={Journal of Computational Dynamics},
  volume={3},
  number={1},
  pages={51--79},
  year={2016},
  publisher={American Institute of Mathematical Sciences}
}

@article{kord18a,
  title={On convergence of extended dynamic mode decomposition to the Koopman operator},
  author={Korda, M. and Mezi{\'c}, I.},
  journal={Journal of Nonlinear Science},
  volume={28},
  number={2},
  pages={687--710},
  year={2018},
  publisher={Springer}
}

@article{brun16a,
  title={Koopman invariant subspaces and finite linear representations of nonlinear dynamical systems for control},
  author={Brunton, S. L. and Brunton, B. W. and Proctor, J. L. and Kutz, J. N.},
  journal={PloS one},
  volume={11},
  number={2},
  pages={e0150171},
  year={2016},
  publisher={Public Library of Science San Francisco, CA USA}
}

@INPROCEEDINGS {Take81,
	AUTHOR = "F. Takens",
	TITLE = "Detecting Strange Attractors in Fluid Turbulence",
	BOOKTITLE = "Symposium on Dynamical Systems and Turbulence",
	EDITOR = "D. A. Rand and L. S. Young",
	JOURNAL = "Lect. Notes Math.",
	PUBLISHER = "Springer-Verlag",
	VOLUME = 898,
	PAGES = 366,
	ADDRESS = "Berlin",
	YEAR = 1981}

@ARTICLE {Pack80,
	TITLE = "Geometry from a Time Series",
	AUTHOR = "N. H. Packard and J. P. Crutchfield and J. D. Farmer and R. S. Shaw",
	JOURNAL = "Phys. Rev. Let.",
	VOLUME = 45,
	PAGES = 712,
	YEAR = 1980}

@article{klus20a,
  title={Data-driven approximation of the {K}oopman generator: Model reduction, system identification, and control},
  author={Klus, S. and N{\"u}ske, F. and Peitz, S. and Niemann, J.-H. and Clementi, C. and Sch{\"u}tte, C.},
  journal={Physica D: Nonlinear Phenomena},
  volume={406},
  pages={132416},
  year={2020},
  publisher={Elsevier}
}

@ARTICLE {Crut87a,
	TITLE = "Equations of Motion from a Data Series",
	AUTHOR = "J. P. Crutchfield and B. S. McNamara",
	JOURNAL = "Complex Systems",
	VOLUME = 1,
	PAGES = "417 -- 452",
	YEAR = 1987}

@article{brun16b,
  title={Discovering governing equations from data by sparse identification of nonlinear dynamical systems},
  author={Brunton, S. L. and Proctor, J. L. and Kutz, J. N.},
  journal={Proceedings of the national academy of sciences},
  volume={113},
  number={15},
  pages={3932--3937},
  year={2016},
  publisher={National Acad Sciences}
}

@article{sugi12a,
  title={Detecting causality in complex ecosystems},
  author={Sugihara, G. and May, R. and Ye, H. and Hsieh, C.-h. and Deyle, E. and Fogarty, M. and Munch, S.},
  journal={Science},
  volume={338},
  number={6106},
  pages={496--500},
  year={2012},
  publisher={American Association for the Advancement of Science}
}

@BOOK {Coll80a,
	TITLE = "Maps of the Unit Interval as Dynamical Systems",
	AUTHOR = "P. Collet and J.-P. Eckmann",
	PUBLISHER = "Birkhauser",
	ADDRESS = "Berlin",
	YEAR = 1980}

@ARTICLE {Miln77,
	TITLE = "On Iterated Maps of the Interval",
	AUTHOR = "J. Milnor and W. Thurston",
	JOURNAL = "Springer Lecture Notes",
	VOLUME = 1342,
	PAGES = "465--563",
	YEAR = 1988}

@book{meiss07a,
  title={Differential dynamical systems},
  author={Meiss, James D},
  year={2007},
  publisher={SIAM}
}

@article{rung19a,
  title={Inferring causation from time series in {E}arth system sciences},
  author={Runge, J. and others},
  journal={Nature communications},
  volume={10},
  number={1},
  pages={1--13},
  year={2019},
  publisher={Nature Publishing Group}
}

@article{arba17a,
  title={Ergodic theory, dynamic mode decomposition, and computation of spectral properties of the {K}oopman operator},
  author={Arbabi, H. and Mezic, I.},
  journal={SIAM Journal on Applied Dynamical Systems},
  volume={16},
  number={4},
  pages={2096--2126},
  year={2017},
  publisher={SIAM}
}

@article{ghil20a,
  title={The physics of climate variability and climate change},
  author={Ghil, Michael and Lucarini, Valerio},
  journal={Reviews of Modern Physics},
  volume={92},
  number={3},
  pages={035002},
  year={2020},
  publisher={APS}
}

\appendix

\section{Random Fourier Features and RKHS} 
\label{app:RFFs}

It is instructive to briefly go through the relation between random Fourier features and kernels. In particular, tensor products of Gaussian-distributed random Fourier features approximate tensor product Gaussian kernels. We largely follow the exposition in Reference \cite{Gund19a}. 

Let $\mathcal{V}$ be a Hilbert space. A function $\varphi:\M \rightarrow \mathcal{V}$ is called a \emph{feature map}.  A feature map naturally defines a kernel, and hence a RKHS via 
\begin{align*}
    \mathsf{K}(\omega, \omega') = \langle \varphi(\omega), \varphi(\omega') \rangle_{\mathcal{V}}
    ~,
\end{align*}
Given a kernel function $\mathsf{K}$ over $\Omega$, there may be many feature maps which represent the kernel: $\mathsf{K}(\omega, \omega') = \langle \varphi(\omega), \varphi(\omega') \rangle_{\mathcal{V}}$.  The map $\varphi(\omega) = k_\omega$ is one example, but in general there will be many more. 


Random Fourier features are random functions $\Psi: \Omega \rightarrow \Reals^M$ that approximate the kernel as
\begin{align*}
    \mathsf{K}(\omega, \omega') = \langle \varphi(\omega), \varphi(\omega') \rangle_{\mathcal{V}} \approx \Psi(\omega)^\top \Psi(\omega')
    ~.
\end{align*}
From Bochner's theorem~\cite{rudi94}, the Fourier transform of a non-negative measure gives a continuous, positive definite, and shift-invariant kernel $\mathsf{K}(\omega, \omega') = \mathsf{K}(\omega - \omega')$, 
\begin{align*}
    \mathsf{K}(\omega - \omega') = \int p(\phi) e^{i \phi^\top(\omega - \omega')} d\phi
    ~.
\end{align*}

Let $\psi(\omega) = e^{i \phi^\top \omega}$ and normalize the non-negative measure $p(\phi)$ so that it is a probability measure. The Fourier transform becomes an expectation over $\phi$ and we take $M$ IID samples $\{\phi_j\}_{j=1}^M$ for a Monte Carlo approximation, 
\begin{align*}
    \mathsf{K}(\omega - \omega') &= \int p(\phi) e^{i\phi^\top(\omega - \omega')} d\phi \\
    &= \mathbb{E}_\phi \bigl[e^{i \phi^\top(\omega - \omega')} \bigr] \\
    &\approx \frac{1}{R} \sum\limits_{j=1}^M e^{i \phi_j^\top(\omega - \omega')} \\
    &= \begin{bmatrix}
        \frac{1}{\sqrt{M}} e^{i \phi_1^\top \omega} \\
        \frac{1}{\sqrt{M}} e^{i \phi_2^\top \omega} \\
        \vdots\\
        \frac{1}{\sqrt{M}} e^{i \phi_M^\top \omega} 
    \end{bmatrix}^\top
    \begin{bmatrix}
        \frac{1}{\sqrt{M}} e^{-i \phi_1^\top \omega} \\
        \frac{1}{\sqrt{M}} e^{-i \phi_2^\top \omega} \\
        \vdots\\
        \frac{1}{\sqrt{M}} e^{-i \phi_M^\top \omega} 
    \end{bmatrix}\\
    &:= \Psi(\omega)\Psi(\omega')^*
    ~.
\end{align*}

Now consider a system with $L$ components, $\omega = [\omega_1 \; \omega_2 \; \cdots \; \omega_L]^\top$. Let $\psi_i(\omega_i) := e^{i \phi_i^\top \omega_i}$ for $i \in \{1, \ldots, L\}$. Taking the tensor product of the component random Fourier feature functions, we get 
\begin{align*}
    [\psi_1 \otimes \psi_2 \otimes \cdots \otimes \psi_L](\omega) &= \prod\limits_{i=1}^L \psi_i(\omega_i) \\ 
    &= e^{i\phi_1^\top \omega_1}e^{i \phi_2^\top \omega_2} \cdots e^{i\phi_L^\top \omega_L} \\ 
    &= e^{i(\phi_1^\top \omega_1 + \phi_2^\top \omega_2 + \cdots + \phi_L^\top \omega_L)} \\ 
    &= e^{i \phi^\top \omega} \\ 
    &= \psi(\omega)
    ~,
\end{align*}
where $\phi = [\phi_1 \; \phi_2 \; \cdots \; \phi_L]^\top$.
Therefore, the tensor product of component random Fourier features is also a random Fourier feature. 

We will use random Fourier features as dictionary basis functions.  From the discussion above, this dictionary approximates the kernel associated to a given shift invariant kernel which arsises from a feature map. This will allow us to create random projections of Koopman operators into a given RKHS function space $\F$. This applies for the full space $\F$ as well as component spaces $\F_i$. 

Since we use Koopman approximations in this work to evolve real-valued functions, we utilize real-valued random Fourier feature dictionaries. Together with some trigonometric manipulations and $b \sim \text{Uniform}(0, 2\pi)$, this gives our random Fourier feature dictionary in Eqn. (\ref{eq:RFF}), 
\begin{align*}
    \Psi_{\text{RFF}}(\omega) = [\psi_1(\omega) \; \cdots \psi_M(\omega)]^\top
    ~,
\end{align*}
with $\psi_j(\omega) = \cos(\phi_j^\top \omega + b)$.

We note that Gaussian radial basis functions (RBFs) also have the above property that they can be used to define both the full RKHS space $\F$ as well as the component spaces $\F_i$. Each RBF dictionary function is given by a different choice of center $\omega'$, which can be chosen, e.g., uniformly or randomly. In practice, we have found it easier to do parameter optimizations with Gaussian random Fourier features than with Gaussian RBFs. Both are valid dictionary choices for the causal discovery algorithms we introduce below, and we have found similar performance when parameters are optimized manually. 

\begin{figure*}
    \begin{center}
        \includegraphics[width=1.0 \textwidth]{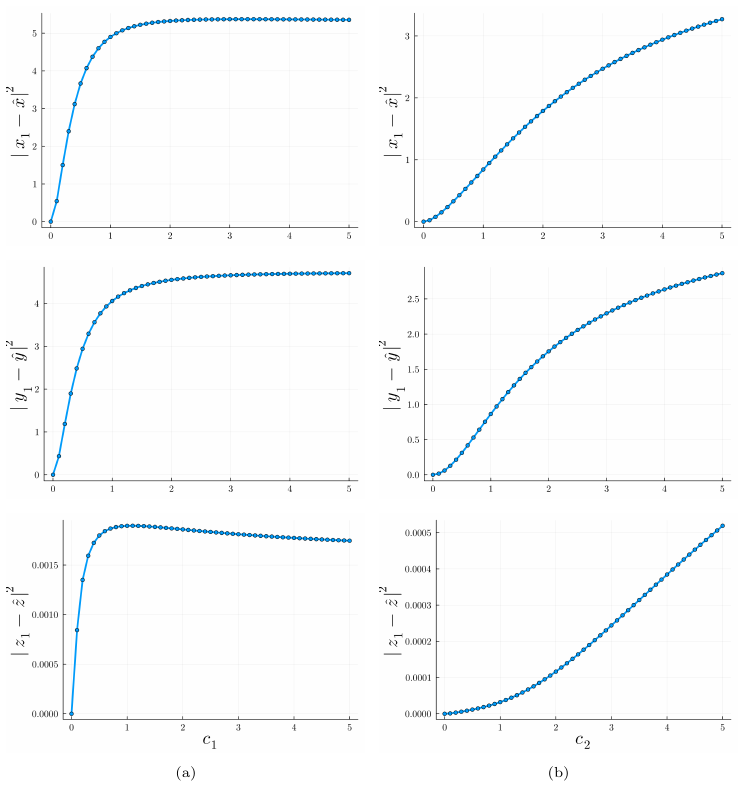}
    \end{center}
    \caption{Some causal phenomenology of the coupled \rossler oscillator system using a counterfactual causality measure. Column (a) shows causality from an asymmetrical system with $\M_2 \causes \M_1$ and $\M_1 \ncauses \M_2$. Column (b) shows similar causality measures for the fully-coupled system $\M_2 \causes \M_1$ and $\M_1 \causes \M_2$.}
    \label{fig:rossler_feedbacks}
\end{figure*}

\section{Feedback Phenomenology of Coupled \rossler System}
\label{app:rossler_feedback}

Here we consider the fully-coupled case of the coupled \rossler system with both $c_1 \neq 0$ and $c_2 \neq 0$. There are feedbacks in this system such that each component has a causal influence on itself at later times, mediated through the other component.  That is, the feedback of $\M_1$ on itself is the combined causal effect of $\M_1 \causes \M_2$ and then $\M_2 \causes \M_1$. 
    
To get some notion of the strength of these feedbacks using counterfactual measures, we will isolate the first piece, $\M_1 \causes \M_2$ (which was not active in Example~\ref{ex:rossler_causal_phenom}, reproduced in Figure~\ref{fig:rossler_feedbacks} (a), because $c_2 = 0$ in that case). Here, we fix $c_1 = 1.0$ and vary $c_2$. For the counterfactual, we no longer compare with an independent $\M_1$ \rossler system (i.e. both $c_1 = 0$ and $c_2 = 0$). Instead, we compare orbits of $\M_1$ in the fully coupled system with orbits from an asymmetrically-coupled system with no feedback, i.e. $c_1 \neq 0$ and $c_2 = 0$. 

This measure as a function of $c_2$ is shown in Figure~\ref{fig:rossler_feedbacks} (b). In the fully coupled case, the dependence on the coupling constant $c_2$ does not saturate in the same ``weak coupling'' regime as in Figure~\ref{fig:rossler_feedbacks} (a). Note also that the magnitude of the causal measures in (b) are systematically lower than in (a), as we would expect given that the causal effects in (b) are mediated through $\M_2$. Unlike the asymmetrically-coupled case, the inclusion of feedbacks prevents the orbits from converging with increased coupling strength---in this case the coupling $c_2$ from $\M_1$ to $\M_2$. This result again emphasizes the need to consider nonlinear systems holistically, even when there is linear coupling and weak nonlinearity. 

\section{Proofs}

\subsection{Proof of Theorem \ref{thm:bounded}}
\label{prf:bounded}
\begin{theorem*}
    Let $\mathcal{F}$ be a RKHS, and $\Koop^t$ satisfy that $\Koop^t f \in \mathcal{F}$ for all $f \in \mathcal{F}$ ($\Koop^t:\mathcal{F} \rightarrow \mathcal{F}$). Then $\Koop^t$ is bounded (continuous).
\end{theorem*}

\begin{proof}
    By the closed graph theorem, $\Koop^t$ is continuous (bounded) if and only if the graph $$\Gamma(\Koop^t):= \{(f, \Koop^t f): f \in \mathcal{F}\} \subset \mathcal{F} \times \mathcal{F}$$
    is closed inside $\mathcal{F} \times \mathcal{F}$.  Equivalently, the graph is closed if and only if $f_i \rightarrow 0$, and $\Koop^t f_i \rightarrow g$, then it must be the case that $g = 0$. Suppose that $f_i \rightarrow 0$ in $\mathcal{F}$.  Then for each $\omega \in \Omega$
    \[
    \Koop^t f_i(\omega) = f_i(\Phi^t(\omega)) = E_{\Phi^t(\omega)}(f_i)
    \]
    where $E_x$ is the evaluation functional at $x \in \Omega$. Notice then that
    \[
    |K^t f_i(\omega)| = |E_{\Phi^t(\omega)}(f_i)| \leq \|E_{\Phi^t(\omega)}\| \|f_i\|_\mathcal{F} \rightarrow 0.
    \]
    Hence, $\Koop^t f_i(\omega) \rightarrow 0$ point-wise. Now by assumption, $\Koop^t f_i \rightarrow g$ in $\mathcal{F}$, and hence, $\Koop^t f_i \rightarrow g$ point-wise as well.  Therefore, $g(\omega) = 0$ for all $\omega$, so that $g=0$, proving the result. 
\end{proof}

\subsection{Proof of Lemma \ref{lemma:isometry}}

\begin{lemma*}
    Let $\M = \M_1 \times \cdots \times \M_L$, $\mathsf{K}_i$, $\mathsf{K}$ and $\F$ be as above.  For each  $X \subset \{1,2, \dots, L\}$,  we have that  $\F_X$  is a  closed subspaces in $\F$.
\end{lemma*}

\begin{proof}
    \label{prf:isometry}
    Clearly $\F_X$ is a subspace of $\F$; it suffices to show it is closed.   Suppose that $(f_n) \subset \F_X$ converges to $f$. Then for all $x \in \M_X$ and $a,b \in \M_{X^c}$, we have
    \[
    \langle f_n, k_{x,a} \rangle = f_n(x,a) = f_n(x,b) = \langle f_n, k_{x,b} \rangle
    \]
    Hence, we have
    \[
    \begin{array}{rcl}
            |f(x,a) - f(x,b)| &=& \left| \langle f, k_{x,a} \rangle -  \langle f, k_{x,b} \rangle \right|\\
            
        &=& \lim\limits_{n \rightarrow \infty} \left| \langle f_n, k_{x,a} \rangle -  \langle f_n, k_{x,b} \rangle \right|\\
    &=&0
    \end{array}
    \]
    Thus, $f \in \F_X$. 
\end{proof}

\subsection{Proof of Proposition \ref{prop:dense span}}

\begin{proposition*}
    Let $X \subset \{1,2, \dots, L\}$, and $\mathsf{K}_i$, $\F_X$ be as above.  Define the kernel function $\mathsf{K}_X: \Omega \times \Omega \rightarrow \mathbb{C}$ by 
    \[
     \mathsf{K}_X(\omega', \omega) := \prod_{j \in X} \mathsf{K}_j(\omega_j', \omega_j)
    \]
    for all $\omega, \omega' \in \Omega$.  Define the functions $k_{X,\omega}: \Omega \rightarrow \mathbb{C}$ by $k_{X,\omega}(\omega') = \mathsf{K}_X(\omega', \omega)$. Then
    \[
    \F_X  = \overline{span_{\omega \in \Omega }} \{k_{X,\omega} \}
    \]
\end{proposition*}

\begin{proof}
    \label{prf:dense span}
    Let $\mathbbm{1}_i$ denote the identity function $\mathbbm{1}_i(\omega_i) = 1$ on $\Omega_i$.  By assumption, $\mathbbm{1}_i \in \F_i$ for each $i$. Note that $k_{X,\omega}$ has the form $k_{X,\omega} = (\otimes_{j \in X} K_{\omega_j}) \otimes (\otimes_{i \in X^c} \mathbbm{1}_i)$.  Indeed, for all $\omega' \in \Omega$, we have 
    \begin{align*}
    K_{X,\omega}(\omega') &= \prod_{j \in X} \mathsf{K}_j(\omega_j', \omega_j)\\ &= \otimes_{j \in X} k_{\omega_j}(\omega_j')\\
    &= \left[ \otimes_{j \in X} k_{\omega_j}(\omega_j')\right]\left[ \otimes_{i \in X^c} \mathbbm{1}_i(\omega_i')\right] \\
    &= (\otimes_{j \in X} K_{\omega_j}) \otimes (\otimes_{i \in X^c} \mathbbm{1}_i) (\omega')
    \end{align*}
    From the form $k_{X,\omega} = (\otimes_{j \in X} K_{\omega_j}) \otimes (\otimes_{i \in X^c} \mathbbm{1}_i)$, we see that $k_{X,\omega}$ is constant on $\Omega_{X^c}$. That is, $k_{X,\omega} \in \mathcal{F}_X$.  Furthermore, we have that for each $j$ that $\overline{\mbox{span}}\{k_{\omega_j}\} = \mathcal{F}_j$.  Using the fact that $\mathbbm{1}_j \in \mathcal{F}_j$ for each $j$, we then get that $\overline{span_{\omega \in \Omega }} \{k_{X,\omega} \} = \F_X$. 
\end{proof}

\subsection{Proof of Theorem \ref{thm:causal_equiv}}

\begin{theorem*}
    Dynamical causal influence and Koopman causal influence are equivalent, 
    \begin{align*}
        \M_C \causes^t_\Koop \M_E \iff \M_C \causes^t \M_E
    \end{align*}
    Correspondingly, 
    \begin{align*}
        \M_C \ncauses^t_\Koop \M_E \iff \M_C \ncauses^t \M_E
    \end{align*}
\end{theorem*}

\begin{proof}
    \label{prf:causal_equiv}
    We will prove $\M_C \ncauses^t_\Koop \M_E \iff \M_C \ncauses^t \M_E$. Let $\omega \in \Omega$. Consider the functions $k_{E,\omega} \in \F_E$ from Proposition \ref{prop:dense span}.   Note that if $\omega = [\omega_C \; \omega_E \; \omega_R]^\top \in \Omega$ then we have
    \begin{equation}
        \label{eq: simple_cause}
              k_{E,\omega}(\Phi^t([\omega_C \; \omega_E \; \omega_R]^\top)) =  \Koop^t k_{E,\omega}([\omega_C \; \omega_E \; \omega_R]^\top)
    \end{equation}
 
    By definition, $\M_C \ncauses^t \M_E$  if and only if $\Phi^t([\omega_C \; \omega_E \; \omega_R]^\top)$ does not depend on $\omega_C$.  From Eqn. \ref{eq: simple_cause}, this happens if and only if the functions $\Koop^t k_{E,\omega}$ do not depend on $\omega_C$.  That happens if and only if $\Koop^t k_{E,\omega} \in \F_{E,R}$. By the linearity of the Koopman operator, we have that for all $f \in \mbox{span}\{k_{E,\omega}\}$, we get that $\M_C \ncauses^t \M_E$ if and only if $\Koop^t f \in \F_{E,R}$. From Theorem \ref{thm:bounded}, $\Koop^t$ is a bounded operator. 
    It follows from the boundedness of the Koopman operator that by passing to limits, we get $\M_C \ncauses^t \M_E$ if and only if  $\Koop^t f \in \F_{E,R}$ for all $f \in \overline{\mbox{span}}\{k_{E,\omega}\}$. However Proposition \ref{prop:dense span} states that $\overline{\mbox{span}}\{k_{E,\omega}\} = \F_E$.  Therefore, we get $\M_C \ncauses^t \M_E$ if and only if  $\Koop^t f \in \F_{E,R}$ for all $f \in \F_E$, which is the definition of $\M_C \ncauses^t_\Koop \M_E$.
\end{proof}

\section{Algorithms and Numerical Study Details}
Here we provide details on our Koopman causality algorithms as they are implemented in our \href{https://github.com/adamrupe/KoopmanCausality}{\texttt{KoopmanCausality.jl}} code written in Julia~\cite{Julia17}. 

Our main algorithm in Sec.~\ref{sec:dmd_causality} A. computes the quantity $\M_C \ddcauses \M_E$, which approximates our measure $\M_C \causes_\Koop^t \M_E$ from Definition~\ref{def:koop_cause}. Recall that this measure defines the causal influence of the `cause' component $\M_C$ on the `effect' component $\M_E$ at the given time-shift $t$. 

The algorithm requires data triples $\bigl(\omega_E(n), \omega_E^t(n), \omega_C(n)\bigr)$ split into \texttt{train} and \texttt{test} sets, where $\omega_E^t := P_E \Phi^t ([\omega_E \omega_C \omega_R]^\top)$ is the time-shift of $\omega_E$ and $n$ is the time index of the data. We collect the \texttt{train} data into data matrices $\mathbf{\M}_E$, $\mathbf{\M}_E^t$, and $\mathbf{\M}_C$, following equations (\ref{eq:data_mat}) and (\ref{eq:data_mat_shift}). 


Because we utilize random Fourier feature dictionaries throughout, we first detail how to create random Fourier feature (RFF) dictionaries $\Psi_{\text{RFF}}$, given in Eqn. (\ref{eq:RFF}). For each of the $M$ cosine functions in our RFF dictionaries, we sample the $N$-dimensional random vector frequencies $\phi_j$ via $N$ IID samples of a zero-mean Gaussian distribution with variance $\sigma$, and the phases $b_j$ are uniformly distributed. Note that we use the same number of random feature functions $M$ in the marginal $\Psi_E$ and joint $\Psi_{E,C}$ dictionaries. One may be tempted to increase the number in the joint dictionary, since the joint vectors are higher dimension, but recall that the function dictionaries are samples of a larger function space, and we want to use the same sampling of the marginal and joint function spaces. 

\begin{algorithm}[t]
\SetAlgoLined
    \Input{Dimensionality $N$, number of features $M$, and variance $\sigma$}
    Initialize \texttt{list} $\Psi_{\text{RFF}}$\;
    \For{each $m$ in $1:M$}{
        create vector $\phi_m$ from $N$ samples of $\text{Normal}(0, \sigma)$\;
        sample $b_m$ from Uniform$[0,1]$ \;
        store random feature function $\psi_m(\omega) = \cos(\phi_m^\top \omega + b_m)$\;
    }
    \Output{$\Psi_{\text{RFF}} = [\psi_1(\omega) \cdots \psi_M(\omega)]$}
    \caption{Create RFF dictionary}
    \label{alg:RFF}
\end{algorithm}

Once we have our chosen dictionaries, we can create the marginal and joint models using DMD. Note that each of these models requires its own dictionary. The dimensionality $N$ of the functions in the marginal model dictionary is the dimensionality of the `effect' component $\text{dim}(\M_E)$, and for the joint model it is the sum $\text{dim}(\M_E) + \text{dim}(\M_C)$. Recall that when using random Fourier features, the joint model dictionary functions take as input the concatenation of the `effect' and `cause' component vectors, shown in Eqn.~(\ref{eq:rff_joint}). 

\begin{algorithm}[!ht]
\SetAlgoLined
  Create the marginal model $\KoopMat^t_{\text{marg}}$\;
  \Input{matrices $\mathbf{\M}_E$, $\mathbf{\M}_E^t$, and dictionary $\Psi_E$}
  \Begin{
      apply $\Psi_E$ to $\mathbf{\M}_E$ to create $\mathbf{\Psi}_E$, (\ref{eq:dict_mat})\;
      concatenate $\mathbf{\Psi}_E$ to $\mathbf{\M}_E$ to create 
      $\begin{bmatrix}
        \mathbf{\Omega}_E \\
        \mathbf{\Psi}_E
    \end{bmatrix}$\;
      multiply $\mathbf{\M}_E^t$ by the pseudoinverse 
      $\begin{bmatrix}
        \mathbf{\Omega}_E \\
        \mathbf{\Psi}_E
    \end{bmatrix}^\dagger$
    \;
      }
    \Output{matrix $\KoopMat^t_{\text{marg}} = 
    \mathbf{\Omega}_E^{t} 
    \begin{bmatrix}
        \mathbf{\Omega}_E \\
        \mathbf{\Psi}_E
    \end{bmatrix}^\dagger$, (\ref{eq:marginal_model})}
    \BlankLine
    
  Create the joint model $\KoopMat^t_{\text{joint}}$\;
  \Input{matrices $\mathbf{\M}_E$, $\mathbf{\M}_E^t$, $\mathbf{\M}_C$ and dict. $\Psi_{E,C}$}
  \Begin{
      apply $\Psi_{E,C}$ to $\begin{bmatrix}\mathbf{\M}_E \\ \mathbf{\M}_C \end{bmatrix}$ to create $\mathbf{\Psi}_{E,C}$, (\ref{eq:dict_mat})\;
      concatenate $\mathbf{\Psi}_E$ to $\begin{bmatrix}\mathbf{\M}_E \\ \mathbf{\M}_C \end{bmatrix}$ to create $\begin{bmatrix}\mathbf{\M}_E \\ \mathbf{\Psi}_{E,C} \end{bmatrix}$\;
       multiply $\mathbf{\M}_E^t$ by the pseudoinverse 
      $\begin{bmatrix}
        \mathbf{\Omega}_E \\
        \mathbf{\Psi}_E
    \end{bmatrix}^\dagger$;
      
      }
  \Output{matrix $\KoopMat^t_{\text{joint}} = 
    \mathbf{\Omega}_E^{t} 
    \begin{bmatrix}
        \mathbf{\Omega}_E \\
        \mathbf{\Psi}_{E,C}
    \end{bmatrix}^\dagger$, (\ref{eq:joint_model})}
\caption{Train models}
\label{alg:train}
\end{algorithm}

With the trained marginal and joint models, we now evaluate their performance over the \texttt{test} set and take the difference in model error to compute our Koopman causality approximation. The key idea, again, is to test whether the joint basis in $\M_E \times \M_C$ better represents the Koopman-shifted test functions $f_E \in \mathcal{F}_E$ than the marginal basis in $\M_E$ only. If so, then $\M_C$ causally influences $\M_E$ over the time interval $t$. 

\begin{algorithm}[!ht]
\SetAlgoLined
    \Input{\texttt{train} data matrices $\mathbf{\M}_E$, $\mathbf{\M}_E^t$, $\mathbf{\M}_C$, \texttt{test} data $\{\omega_E, \omega_E^t, \omega_C\}_{\texttt{test}}$, and dictionaries $\Psi_E$, $\Psi_{E,C}$}
    Create models $\KoopMat^t_{\text{marg}}$ and $\KoopMat^t_{\text{joint}}$, Alg.~\ref{alg:train}\;
    \BlankLine
    Evaluate marginal model error\;
    \For{each $(\omega_E, \omega_E^t)$ in \texttt{test}}{
        generate model prediction $\widetilde{\omega}_E^{t} |_{\text{marg}}$, Eqn.~(\ref{eq:marg_pred})\;
        compute square error $||\widetilde{\omega}_E^{t} |_{\text{marg}} - \omega_E^t||^2$ \;
    }
    marginal error $= \frac{1}{|\texttt{test}|} \sum\limits_\texttt{test} ||\widetilde{\omega}_E^{t} |_{\text{marg}} - \omega_E^t||^2$\;
    \BlankLine
    Evaluate joint model error\;
    \For{each $(\omega_E, \omega_E^t, \omega_C)$ in \texttt{test}}{
        generate model prediction $\widetilde{\omega}_E^{t} |_{\text{joint}}$, Eqn.~(\ref{eq:joint_pred})\;
        compute square error $||\widetilde{\omega}_E^{t} |_{\text{joint}} - \omega_E^t||^2$ \;
    }
    joint error $= \frac{1}{|\texttt{test}|} \sum\limits_\texttt{test} ||\widetilde{\omega}_E^{t} |_{\text{joint}} - \omega_E^t||^2$\;
    \BlankLine
    \Output{$\M_C \ddcauses \M_E = \text{marginal error}-\text{joint error}$}
    \caption{Approximate Koopman Causality}
    \label{alg:koop_cause}
\end{algorithm}

Finally, we detail conditional forecasting. This proceeds similar to our Koopman causality measure, with some slight modification. It uses trained marginal and joint models as given in Algorithm~\ref{alg:train}, but with a \textbf{small time-shift} $t$, typically a single step in discrete time. The main difference is how the model is evaluated over the \texttt{test} set. For conditional forecasting, prior-step predictions are combined with dictionary evaluations of values in the \texttt{test} set to create the next-step predictions. 

\begin{algorithm}[!ht]
\SetAlgoLined
    \Input{\texttt{train} data matrices $\mathbf{\M}_E$, $\mathbf{\M}_E^t$, $\mathbf{\M}_C$, \texttt{test} data $\{\omega_E, \omega_E^t, \omega_C\}_{\texttt{test}}$, and dictionaries $\Psi_{\text{marg}}$, $\Psi_{\text{joint}}$}
    Create models $\KoopMat^t_{\text{marg}}$ and $\KoopMat^t_{\text{joint}}$, Alg.~\ref{alg:train}\;
    set $\widetilde{\omega}_E(1) = \omega_E(1)$\;
    \BlankLine
    Marginal model rollout\;
    \For{each $i, \omega_E$ in enumerate(\texttt{test})}{
        generate model prediction $\widetilde{\omega}_E(i+1) |_{\text{marg}} = \KoopMat^t_{\text{marg}}
    \begin{bmatrix}
        \widetilde{\omega}_E(i) \\
        \Psi_E\bigl(\omega_E(i)\bigr)
    \end{bmatrix}$\;
    }
    \BlankLine
    Joint model rollout\;
    \For{each $i, (\omega_E, \omega_C)$ in enumerate(\texttt{test})}{
        generate model prediction $\widetilde{\omega}_E(i+1) |_{\text{joint}} = \KoopMat^t_{\text{joint}}
    \begin{bmatrix}
        \widetilde{\omega}_E(i) \\
        \Psi_{E,C}\bigl(\omega_E(i), \omega_C(i)\bigr)
    \end{bmatrix}$\;
    }
    \Output{Conditional forecasts $\{\widetilde{\omega}(i)\}_{\text{marg}}$ and $\{\widetilde{\omega}(i)\}_{\text{joint}}$}
    \caption{Conditional Forecasting}
    \label{alg:cond_forec}
\end{algorithm}

\subsection{Random Fourier Feature Dictionary Optimization}
The most impactful hyperparameter when using random Fourier feature dictionaries is the variance $\sigma$ in Algorithm~\ref{alg:RFF}. To find an optimal value of $\sigma$, we utilize a naive brute-force search by choosing a range of values for $\sigma$ and then train and evaluate models for each of these values. Specifically, we use $\{\sigma\} = \{0.0005, 0.001, 0.005, 0.01, 0.05, 0.1, 0.25, 0.5, 1.0, 2.0\}$. Moreover, because the features are sampled randomly for a given value of $\sigma$, we actually train and evaluate multiple ($\eta$) models for each $\sigma$. In this work, we train and evaluate $\eta = 5$ marginal and joint models for each value of $\sigma$ in our hyperparameter optimization loop. After training and evaluating all the marginal and joint models in the optimization loop, we keep those with the minimal marginal error and joint error, respectively. Within the respective model search spaces, those with minimal error best show how the marginal dictionary and joint dictionary function spaces can represent the Koopman evolution of the test functions $f_E \in \mathcal{F}_E$ (which we again take to be the identity $f_E(\omega_E) = \omega_E$). Finally, because the number of random feature functions $M$ in the dictionary does not have much of an effect on the outcome of our algorithm, we pick a fixed value $M = 150$ for all cases. 

\begin{algorithm}[t]
\SetAlgoLined
    \Input{Set of variances $\{\sigma\}$, number of samples $\eta$, number of features $M$}
    Initialize \texttt{marg\_err} = Inf\;
    Initialize \texttt{joint\_err} = Inf\;
    \For{each $\sigma$ in $\{\sigma\}$}{
        \For{\_ in 1:$\eta$}{
            create $\Psi_E$ with $\sigma$\;
            train and evaluate marginal model with $\Psi_E$ to get \texttt{cur\_marg\_err}, see Alg~\ref{alg:koop_cause} \;
            \texttt{marg\_err} = min(\texttt{marg\_err}, \texttt{cur\_marg\_err})\;
            \BlankLine
            create $\Psi_{E,C}$ with $\sigma$\;
            train and evaluate joint model with $\Psi_{E,C}$ to get \texttt{cur\_joint\_err}, see Alg~\ref{alg:koop_cause}\;
            \texttt{joint\_err} = min(\texttt{joint\_err}, \texttt{cur\_joint\_err}\;
        }
    }
    \Output{\texttt{marg\_err} - \texttt{joint\_err}}
    \caption{RFF causality optimization}
    \label{alg:RFF_opt}
\end{algorithm}

\subsection{Numerical study details}
Here we give details of the numerical studies shown throughout. Jupyter notebooks that generate the results and figures can be found on our \href{https://github.com/adamrupe/KoopmanCausality}{\texttt{KoopmanCausality.jl}} GitHub. 

\subsubsection{Figure 1}
The uni-directional coupled \rossler model used for this counterfactual causality analysis has parameters: $\varphi_1 = 0.8$, $\varphi_2 = 0.5$, $a = 0.2$, $b = 0.2$, $d = 5.7$, and $c_2 = 0.0$, with $c_1$ varying over the values shown. The initial condition is $\omega_0 = (1.0, -2.0, 0.0, 0.11, 0.2, 0.1)$, the orbits are evolved over the time interval $(0.0, 10.0)$ with no transient time, and the standard $\delta = 0.01$ integration time step is used. 

\subsubsection{Figure \ref{fig:rossler_causes_per_time}}
The uni-directional coupled \rossler model for this Koopman causality measure over time analysis has parameters: $\varphi_1 = 0.8$, $\varphi_2 = 0.5$, $a = 0.2$, $b = 0.2$, $d = 5.3$, $c_1 = 0.5$, and $c_2 = 0.0$. The initial condition $\omega_0 = \Phi^t(\omega^*)$ is generated by a transient-time evolution of $t=2000$ from $\omega^* = (1.0, -2.0, 0.0, 0.11, 0.2, 0.1)$, and the integration step is $\delta = 0.01$. The \texttt{train} set for each Koopman causality evaluation is 500,000 time steps and each \texttt{test} set is 10,000 time steps. The marginal and joint models both use $M = 120$ random Fourier features with hyperparameter optimization performed over $\{\sigma\} = \{0.5, 0.1, 0.05, 0.01, 0.001\}$.

\subsubsection{Figure \ref{fig:rossler_causalDMD}}
The uni-direction coupled \rossler model used for conditional forecasting is the same as above in Figure~\ref{fig:rossler_causes_per_time}. The \texttt{train} set is 200,000 time steps and the \texttt{test} set is 2,000 time steps. The DMD Koopman models use $M=120$ random Fourier features. Because these experiments are much quicker to run than, e.g., those in Figure~\ref{fig:rossler_causes_per_time}, hyperparameter search for the optimal $\sigma$ was performed manually, using visual inspection of the forecasts shown in Figure~\ref{fig:rossler_causalDMD}. The marginal model in Figure~\ref{fig:rossler_causalDMD} (a) uses $\sigma = 1.0$, the joint model in (b) uses $\sigma = 0.1$, the marginal model in (c) uses $\sigma = 0.5$, and the joint model in (d) uses $\sigma = 0.5$. 

\subsubsection{\rossler Confounder}
We use the following parameters for both the symmetric and asymmetric \rossler confounder models Section~\ref{sec:confounding}: $\varphi_1 = 0.8$, $\varphi_2 = 0.5$, $\varphi_3 = 0.5$, $a = 0.2$, $b = 0.2$, and $d = 5.3$. For the symmetric case, $c_2 = c_3 = 0.5$, and for the asymmetric case $c_2 = 0.5$ and $c_3 = 0.3$. In both cases, the initial condition is $\omega_0 = \Phi^{1000}(\omega^*)$ with $\omega^* = (1.0, -2.0, 0.0, 0.11, 0.2, 0.1, -0.11, 0.3, 0.0)$. The \texttt{train} sets are 800,000 time steps and the \texttt{test} sets are 20,000 time steps. The Koopman causality models all use the standard random Fourier feature dictionaries with hyperparameter optimization given above in Appendix D 1. 

\subsubsection{Lorenz 96}
For all Lorenz 96 model experiments, we use a moderate forcing parameter $F=4.0$ (the models are not chaotic for this forcing value) and random initial conditions where each of the $N=101$ oscillators is randomly from Uniform$[0,1]$. 

The perturbation propagation in Figure~\ref{fig:L96_pert} is created using two initial conditions that only differ at the middle oscillator $n=51$ by a unit perturbation $\delta = 1.0$. 

For the cumulative causal effects shown in Figure~\ref{fig:L96_cause} and Figure~\ref{fig:L96_instant}, we use initial conditions $\omega_0 = \Phi^{2000}(\omega^*)$ with random $\omega^*$. The \texttt{train} sets are 500,000 time steps and the \texttt{test} sets are 10,000 time steps. The Koopman causality models use $M=120$ random Fourier features with hyperparameter optimization performed over $\{\sigma\} = \{2.0, 1.0, 0.5, 0.1, 0.05, 0.01, 0.005, 0.001\}$. 

\section{Method comparison details}
\label{app:comparison}
For our comparison, we use the coupled \rossler model in Eqn.~(\ref{eq:rossler}). To achieve uni-directional coupling with $\M_2 \causes \M_1$ but $\M_1 \ncauses \M_2$, we set $c_1 = 0.5$ and $c_2 = 0.0$. The other model parameters used are: $\varphi_1 = 0.8$, $\varphi_2 = 0.5$, $a = 0.2$, $b = 0.2$, and nonlinearity $d = 5.3$. The integration time step is $0.01$, the \texttt{train} set is 500,000 time steps long, and each of the 100 \texttt{test} sets is 20,000 time steps long. To randomly sample initial conditions on the attractor for the different \texttt{test} runs, we start with a fixed $\omega^* = (1.0, -2.0, 0.0, 0.11, 0.2, 0.1)$ and then randomly sample a transient evolution time $\tau$ so that each random initial condition is $\omega_0 = \Phi^\tau(\omega^*)$.

\subsection{Koopman causality}
To evaluate our Koopman causality measure $\M_C \ddcauses \M_E$ in this comparison we use a time shift of $t = 1000$ time steps. For each \texttt{test} run, we utilize random Fourier feature dictionaries with the hyperparameter optimization scheme described above. If $\M_C \ddcauses \M_E > 0.1$ we determine that our method has identified a causal relation, and if $\M_C \ddcauses \M_E \leq 0.0$ our method has identified no causal relation. Although we did not encounter any in our experiments, we would consider $\M_C \ddcauses \M_E$ to be indeterminate for values in the range $(0.0, 0.1]$.

\subsection{Granger causality}

Granger causality is a statistical hypothesis test introduced by Clive Granger to infer directional relationships between time series based on predictability \cite{gran69}. Abstractly, a variable \( Y \) is said to \emph{Granger-cause} a variable \( X \) if the past values of \( Y \) improve the prediction of \( X \) beyond what is possible using only the past values of \( X \). Practically, this is implemented under the assumption of a linear relationship across time with second-order stationarity. The test is typically carried out by fitting a vector autoregressive (VAR) model and testing whether lagged values of \( Y \) provide statistically significant explanatory power for \( X \), using a Wald-type F-test.

The Granger causality test is implemented using the \texttt{VAR} model from the \texttt{statsmodels} Python package \cite{seabold2010statsmodels}. Given two multivariate time series \( X(t), Y(t) \in \mathbb{R}^d \), the method stacks them into a single matrix of shape \( (T, 2d) \) and fits a VAR model with a user-defined maximum lag \( \ell_{\max} \). Causality is assessed in both directions by applying Wald tests: (i) testing whether lagged values of \( Y \) contribute to predicting \( X \) and (ii) testing the reverse. Each test yields a p-value, and causality is inferred when the p-value is below a chosen significance level \( \alpha \) (default \( 0.05 \)). For the \rossler ensemble datasets, the test is repeated independently across realizations to count how often each direction of causality is statistically significant, providing a frequency-based summary of causal directionality.  Testing across hyperparameters \(\ell_{\max} \in \{1, \dots, 20\}\) yielded no difference in results.  

\subsection{Convergent cross mapping}
Convergent Cross Mapping (CCM) is a nonlinear causality inference technique introduced by Sugihara et al.\ to study causal relationships in complex ecological systems \cite{sugihara2012detecting}. The method is grounded in \emph{Takens' embedding theorem}, which states that it is possible to reconstruct the state space of a deterministic dynamical system using time-delay embeddings of a single observed variable. Specifically, given a scalar time series \( x(t) \), one can construct a delay-coordinate map
\[
M_x(t) = \left[ x(t), x(t - \tau), \dots, x(t - (E-1)\tau) \right]
\]
which faithfully represents the underlying attractor, provided that the embedding dimension \( E \) is sufficiently large. The manifold reconstructed in this way, known as the \emph{shadow manifold}, is diffeomorphic to the original system’s attractor, meaning it preserves the essential topological and dynamical structure. Consequently, if a variable \( X \) causally influences a variable \( Y \), then the state of \( X \) can be inferred from the shadow manifold \( M_Y \). This \emph{cross mapping} property enables one to test for causality by attempting to reconstruct \( X \) from \( Y \), and vice versa, and checking whether predictive skill increases with the library size. If \( X \rightarrow Y \), then \( M_Y \) will contain information about \( X \), but not necessarily the reverse.

Given two multivariate time series \( X(t), Y(t) \in \mathbb{R}^d \), the CCM algorithm begins by constructing a shadow manifold \( M_Y(t) \) from delay embeddings of \( Y \) using embedding dimension \( E \), time lag \( \tau \), and a total of \( L \) time steps (the library size). Each embedded point \( M_Y(t) \in \mathbb{R}^{E \cdot d} \) is formed by stacking delayed vectors \( [Y(t), Y(t - \tau), \dots, Y(t - (E-1)\tau)] \). For each such point, the algorithm finds its \( E+1 \) nearest neighbors in the manifold, and uses the corresponding time indices to extract \( X(t) \) values. These are combined via a distance-weighted average to estimate \( \hat{X}(t) \), using weights that decay exponentially with distance. The predictive skill is then measured by the \emph{Pearson correlation coefficient} between \( X(t) \) and \( \hat{X}(t) \), averaged across spatial dimensions. A high correlation indicates that information about \( X \) is embedded in \( Y \), supporting a causal interpretation of \( X \rightarrow Y \). To determine directionality, the procedure is repeated for \( Y \) reconstructed from \( M_X \).

For the directional \rossler System, we applied CCM across the range of hyperparameters \(\tau = 1\), \(E \in \{2, \dots 10\}\) and \(L \in \{1,\dots ,4000\}\). If we detect a Pearson correlation coefficient above our chosen threshold of 0.1, we say that causality is detected from \(\Omega_1 \rightarrow \Omega_2\).  Causality from \(\Omega_2 \rightarrow \Omega_1 \) is then checked.  This  is repeated for each ensemble member.  We note that for nearly all ensemble members and all hyperparameter choices, the correlation coefficient was very close to 1 for both \(\Omega_1 \rightarrow \Omega_2\) and \(\Omega_2 \rightarrow \Omega_1\), so our choice of threshold was not very important in the overall results. 

\end{document}